\documentclass[11pt,english]{amsart}%draft,
\usepackage{amsfonts,amsmath,amsthm,amscd,amssymb,latexsym,amsbsy}
\usepackage[dvips]{graphicx}

\DeclareMathOperator*{\csh}{cosh} 
\DeclareMathOperator*{\tnh}{tanh} \DeclareMathOperator*{\tng}{tan}
\DeclareMathOperator{\supp}{supp} \DeclareMathOperator{\sgn}{sgn}

\begin{document}

\newcommand{\norm}[1]{\left\Vert#1\right\Vert}
\newcommand{\abs}[1]{\left\vert#1\right\vert}

\theoremstyle{plain}
\newtheorem{g_theorem}{Theorem}[section]
\newtheorem{g_lemma}{Lemma}[section]
\newtheorem{g_definition}{Definition}[section]
\newtheorem{g_corollary}{Corollary}[section]
\newtheorem*{g_hyp}{Hypothesis}

\theoremstyle{definition}
\newtheorem{g_remark}{Remark}[section]

\numberwithin{equation}{section}

%\date{({\it Presented by M. M. Malamud})}

\title[Solvable models for the Schr\"{o}dinger operators]
{Solvable models for the Schr\"{o}dinger operators with $\delta'$-like potentials}

\author{Yuriy D. Golovaty, Stepan S. Man'ko}%
\address{ Department of Mechanics and Mathematics, Ivan Franko National University of L'viv, 1 Universytetska str., 79000 L'viv, Ukraine }

\begin{abstract}
We turn back to the well known problem of interpretation of the Schr\"{o}dinger operator with the pseudopotential $\alpha\delta'(x)$, where $\delta(x)$ is the Dirac function and $\alpha\in \mathbb{R}$. We show that the problem in its conventional formulation contains  hidden parameters and  the choice of the proper selfadjoint operator is ambiguously determined.
We study the asymptotic behavior of spectra and eigenvectors of the Hamiltonians with increa\-sing smooth potentials  perturbed by  short-range potentials  $\alpha\varepsilon^{-2}\Psi(\varepsilon^{-1}x)$. Appropriate solvable models  are constructed and the corresponding approximation theorems are proved.
We introduce the concepts of the resonance set $\Sigma_\Psi$ and the coupling function
$\theta_\Psi \colon\Sigma_\Psi \to \mathbb{R}$,
which are spectral characteristics of  the potential
$\Psi$. The selfadjoint operators in the solvable models are determined by means of the resonance set  and the coupling function.
\end{abstract}

\subjclass{34L40, 34B09, 81Q10}

\keywords{Schr\"{o}dinger operator, point interaction, pseudopotential, $\delta'$-potential,
solvable model, Krein space, selfadjoint extension, asymptotic expansion}

%\umblabel

\maketitle

\section*{Introduction}

In this paper we investigate the one-dimensional Schr\"{o}dinger operators with the coefficients having singular support on a disjoint set of points. Such ope\-rators are named \textit{operators with singular point interactions} and are closely related to exactly solvable problems in quantum mechanics, atomic physics and acoustics.
Such point interaction models are \textit{solvable} in the sense that the resolvents of corresponding operators can be given explicitly. As a consequence the spectrum and scattering quantities can also be determined explicitly.

The main question here is how to interpret the differential operators with coefficients equal to the generalized functions since  space of distributions  $\mathcal{D}'(\mathbb{R}^n)$ is not
an algebra with respect to the ``pointwise'' multiplication. Generally the problem of choosing a single Hamiltonian corresponding to a particular singular point interaction is complicated enough.
The theory of one-dimensional Schr\"{o}dinger operators with singular potentials has been summarized in research monographs \cite{Albeverio2edition, AlbeverioKurasov}.

According to the postulates of quantum mechanics every physically measurable quantity is specified by a \textit{selfadjoint} operator in a Hilbert space. Starting with  the pioneering work of F. Berezin and L. Faddeev \cite{BerezinFadeev}, the singular point interactions
 have been studied with the help of the selfadjoint extensions theory for symmetric
operators.
As a rule, the symmetric operators in physical models are generated by differential operators.
This leads to the question how to describe all selfadjoint extensions in terms of the boundary value conditions.
In the case of ordinary differential operators the problem has been solved by M.~Krein \cite{AkhiezerGlazman, Krein}.

For partial differential operators the problem becomes more complicated by the infinite deficiency indices of symmetric operators. In that case new approach to the theory of selfadjoint extensions is based on
the concept of boundary triples \cite{Gorbachuky, GorbachukyKochubej} and the corresponding Weyl
functions. In \cite{DerkachMalamud1, DerkachMalamud2} the concept of Weyl function
was generalized to an arbitrary symmetric operator  with infinite deficiency indices.
It is well known that the Weyl function can be used to investigate the spectrum of a selfadjoint extension
and the scattering matrix \cite{BrascheMalamudNeidhardt2002}.
It is worth to note that the concept of boundary triples is useful in the investigation of the point interaction mo\-dels \cite{GoloschapovaOridoroga, KochubeyUMZh1989,
KochubeySibMatZh1991, KochubeyMatNotes1979}.

The singular finite rank perturbations have been studied in papers \cite{AlbeverioKoshmanenko1999, KuzhelNizhnik,
KoshmanenkoUMZh1991, Nizhnik2001} and monograph \cite{KoshmanenkoBook}. The Schr\"{o}dinger operators with different singular potentials have been investigated in \cite{HrynivMykytiuk2001, MikhailetsMolyboga, NizhFAA2003,
ShkalikovSavchMatNotes2006, ShkalikovSavchukMatNotes1999,
ShkalikovSavchukTMMO2003}.

In some instances the set of all selfadjoint extensions of a symmetric operators is rich enough.
Therefore the harder question comes: how to choose an extension that is best suited to our physical model.
For some models the proper operator can not be chosen within the selfadjoint extensions theory, because
the models contain hidden parameters. If we replace the singular potential with a sequence of short-range smooth potentials, then the operator obtained in the zero-range limit, as often happens, can depend on the type of regularization, that is to say, the operator is governed by the profile of squeezed potentials. This profile is a hidden parameter and plays a crucial role in the choice of a selfadjoint extension corresponding to the physical model under consideration.

Differential operators with singular coefficients can be studied within the framework of the contemporary theories such as
the Colombeau algebras \cite{Colombeau} or the  Egorov  generalized  function algebra \cite{Egorov}.
These theories allow us not only to multiply the generalized functions, but they also contain  rich sets of
``$\delta$-functions with  fixed profiles''  as  different elements of the algebra.
This approach to study of  solvable models has been used in \cite{Antonevich1, Antonevich2}.

\section{Problem Statement and Main Results}
\subsection{Problem statement}

In this work we turn back to the well-known problem how to define the one-dimensional Schr\"{o}dinger operator with the $\delta'$-potential. New attempt will be made to answer this question.

Let us consider the formal Hamiltonian
$$H_\alpha=-\frac{d^2}{dx^2}+U(x)+\alpha\delta'(x),\qquad x\in\mathbb R,$$
where $U$ is a real smooth function, $\delta'$ is the first derivative of the Dirac delta function, and $\alpha$ is a real constant. If product $\delta'(x)y(x)$ is treated as  $y(0)\delta'(x)-y'(0)\delta(x)$, then for $\alpha\neq 0$ there exist no solutions in $\mathcal{D}'(\mathbb{R})$ to equation $H_\alpha y=\lambda y$, except for a trivial one.
In order for operator  $H_\alpha$ to be assigned a meaning, firstly we must describe all selfadjoint extensions in
 $L_2(\mathbb{R})$  of the symmetric operator
\begin{equation*}
    L=-\frac{d^2}{dx^2}+U(x),\qquad \mathcal{D}(L)=\{f\in C_0^\infty(\mathbb{R})\colon f(0)=f'(0)=0\}.
\end{equation*}
Next, one of the extensions must be prescribed for the Schr\"{o}dinger operators with potential
$U(x)+\alpha\delta'(x)$.
At present, in both physics- and mathematics-oriented literature there is no a consensus regarding
the pro\-per choice of such operator, since  $L$ has a rich set of selfadjoint extensions (see Sec.~\ref{SectionHistory} for the historical remarks).

Let us consider more realistic model, namely the family of Hamiltonians with smooth short-range potentials, which
approximate the singular potential $\alpha\delta'(x)$. Let
$\mathcal{H}_{\varepsilon}(\alpha,\Psi)$ be the closure in  $L_2(\mathbb{R})$ of the essentially selfadjoint operator \cite[p. 50]{BS}
\begin{equation*}
    H_{\varepsilon}(\alpha, \Psi)=-\frac{d^2}{dx^2}+U(x)+\frac{\alpha}{\varepsilon^{2}}\Psi(\varepsilon^{-1}x), \quad
\mathcal{D}(H_{\varepsilon}(\alpha, \Psi))=C_0^\infty(\mathbb{R}).
\end{equation*}
Here $\varepsilon$ is a small positive parameter. We call $\Psi\in C^\infty_0(\mathbb{R})$ the \emph{profile} of a local perturbation, and  $\alpha$ the \emph{coupling constant}. It is clear that for some profiles  sequence
${\varepsilon^{-2}}\Psi(\varepsilon^{-1}x)$ converges to $\delta'(x)$ in $\mathcal{D}'(\mathbb{R})$.
We will denote by $\mathcal{P}$ the set of real functions $\Psi\in C^\infty_0(\mathbb{R})$ such that $\supp{\Psi}=[-1,1]$. Let $\mathcal{E}(L)$ denote the set of all selfadjoint extensions of $L$.
Suppose also that  potential $U$ increases as $|x|\to \infty$. It follows that the spectrum of $\mathcal{H}_{\varepsilon}(\alpha,\Psi)$ is discrete.

Our main goal is to construct the map $\mathbb R\times
\mathcal{P}\longrightarrow \mathcal{E}(L)$ that assigns a selfadjoint extension  $\mathcal{H}(\alpha,\Psi)$ of operator $L$ to each pair $(\alpha,\Psi)$. The choice of an operator is determined by the proximity
of the energy levels and the pure states for the Hamiltonians with smooth and singular potentials respectively.
We explore the asymptotic behavior of eigenvalues and eigenfunctions as $\varepsilon\to 0$. The asymptotics provides the limit operator  $\mathcal{H}(\alpha,\Psi)$.

\subsection{Structure of Paper}

The paper is organized as follows. In Sec.~\ref{SectionSpectrumStructure}, we derive  qualitative properties of the discrete spectrum of
perturbed operators  $\mathcal{H}_\varepsilon(\alpha,\Psi)$. We show that all eigenvalues are  continuous functions on  $\varepsilon$ that are bounded from above. Generally speaking, the spectrum of this family is not bounded from below. For some profiles $\Psi$ and $\alpha\neq 0$ there exists a finite number of eigenvalues which go to  $-\infty$ as $\varepsilon\to 0$. By the \emph{bounded spectrum} of $\mathcal{H}_\varepsilon(\alpha,\Psi)$ we mean the set of such eigenvalues  that remain bounded as $\varepsilon\to 0$.

The leading terms of asymptotic expansions for the eigenvalues of the bounded spectrum  and the limit operators $\mathcal{H}(\alpha,\Psi)$ are formally constructed in Sec.~\ref{SectionMainTerms}. We introduce two spectral characteristics of profile $\Psi$, namely the \emph{resonance set} $\Sigma_\Psi$ that is the spectrum of the Sturm--Liouville problem $-w''+\alpha \Psi w=0$, $w'(-1)=0$, $w'(1)=0$ on interval $(-1,1)$ with respect to spectral parameter $\alpha$,  and the \emph{coupling function}  $\theta_\Psi
\colon\Sigma_\Psi \to \mathbb{R}$. We set $\theta_\Psi (\alpha)=w_\alpha(1)\bigl(w_\alpha(-1)\bigr)^{-1}$, where
$w_\alpha$ is an eigenfunction corresponding to eigenvalue $\alpha\in \Sigma_\Psi$.
In the case, when the coupling constant doesn't belong to the resonant set, $\mathcal{H}(\alpha,\Psi)$ is just the direct sum  of the  Schr\"odinger operators with potential $U$
on half-axes subject to the Dirichlet boundary condition at the origin.
In the resonant case, when $\alpha\in \Sigma_\Psi$, operator $\mathcal{H}(\alpha,\Psi)$ is the selfadjoint extension of $L$ that corresponds to the coupling conditions
$ f(+0)=\theta_\Psi (\alpha)f(-0)$, $\theta_\Psi
(\alpha)f'(+0)=f'(-0)$ at the origin.

In Sec.~\ref{SectionHistory} we provide a detail discussion of the previous results dealing with
the $\delta'$-interactions and the $\delta'$-potentials.
We give two examples of exactly solvable models with  piecewise constant $\delta'$-like potentials for which
the resonance set and the coupling function are  explicitly computed.
We solve the scattering problem for this $\delta'$-like potential. The limit value of the transmission coefficient
can be also explicitly computed in terms of  coupling function $\theta_\Psi$.

The remainders of asymptotics for the eigenvalues and eigenfunctions of $\mathcal{H}_\varepsilon(\alpha,\Psi)$ are constructed in Sec.~\ref{SectionCorrector}, because we are in need of more precise asymptotics, in order to proof the approximation theorems.

Section~\ref{SectionMainResult} contains the main results of this paper. We prove in this section that the eigenvalues of $\mathcal{H}_\varepsilon(\alpha,\Psi)$ that belong to the bounded spectrum converge to  eigenvalues of $\mathcal{H}(\alpha,\Psi)$. The convergence in $L_2(\mathbb R)$ of the corresponding eigenfunctions is also established.

\section{Spectrum of  $\mathcal{H}_\varepsilon(\alpha,\Psi)$ and  Auxiliary Results}\label{SectionSpectrumStructure}

Firstly we describe the set  $\mathcal{E}(L)$ of all selfadjoint extensions of  minimal operator $L$. The adjoint operator $L^*=-\frac{d^2}{dx^2}+U(x)$ is defined on the domain
$$
\mathcal{D}(L^*)=\{v\in W_2^2(\mathbb{R}\setminus 0)\colon -v''+Uv\in L_2(\mathbb{R})\}.
$$
\begin{g_lemma}\label{AllExtensionsLemma}
Every element of  $\mathcal{E}(L)$ coincides with operator $L^*$, restricted to the set of functions, satisfying the boundary conditions at the origin of one of the types
\begin{equation}\label{SeparatedCond}
      h_1^- v'(-0)=h_2^- v(-0),\qquad
      h_1^+ v'(+0)=h_2^+ v(+0)
\end{equation}
with the parameters $(h_1^\pm, h_2^\pm)$ from the projective space $\mathbb{P}^1$;
\begin{equation}\label{ConnectedCond}
\begin{pmatrix} v(+0) \\ v'(+0) \end{pmatrix}
=C\begin{pmatrix} v(-0) \\ v'(-0) \end{pmatrix},\qquad
C=e^{i\varphi}\begin{pmatrix} c_{11} & c_{12} \\ c_{21} & c_{22} \end{pmatrix},
\end{equation}
where $\varphi\in
[-\frac{\pi}{2},\frac{\pi}{2}]$, $c_{kl}\in\mathbb{R}$ and $c_{11}c_{22}-c_{12}c_{21}=1$.
\end{g_lemma}
This lemma has been proved in \cite{ChernoffHughesJFunctAnal93, SebaCzechJPhys86} for $U=0$.
Obviously, the smooth potential does not affect on the boundary conditions at the origin.
Selfadjoint operators described by  boundary conditions \eqref{SeparatedCond} will be named
\textit{separated extensions}. These operators are equal to the orthogonal sum of
two selfadjoint operators defined on the half-axes.
Selfadjoint ope\-ra\-tors of the second type will be named  \textit{connected extensions},
because these conditions connect the boundary values of the function on
the left and right half-axes.

Let us write $\Psi_\varepsilon(x)=\varepsilon^{-2}\Psi(\varepsilon^{-1}x)$ and set $m_k(f)=\int_{\mathbb{R}}\xi^k f(\xi)\,d\xi$.
\begin{g_lemma}\label{LemDelFunc}
Assume $\Psi\in \mathcal{P}$ and  $c$ is a nonzero constant. Sequence $\Psi_\varepsilon$ converges to $c\delta'(x)$ as $\varepsilon\to 0$ in the sense of distributions  if and only if
\begin{equation}\label{CondDeltPrime}
    m_0(\Psi)=0, \quad m_1(\Psi)\neq 0.
\end{equation}
In addition, $c=-m_1(\Psi)$.
\end{g_lemma}
\begin{proof}
Given $\varphi\in C_0^\infty(\mathbb{R})$, we have
\begin{multline*}
\langle\Psi_\varepsilon,\varphi\rangle=
\varepsilon^{-2}\int\limits_{-\varepsilon}^{\varepsilon}
\Psi(\varepsilon^{-1}x)\varphi(x)\,dx=
\varepsilon^{-1}\int\limits_{-1}^1 \Psi(\xi)\varphi(\varepsilon\xi)\,d\xi\\
=\varepsilon^{-1}\int\limits_{-1}^1
\Psi(\xi)\bigl(\varphi(0)+\varepsilon\varphi'(0)\xi+O(\varepsilon^2)\bigr)\,d\xi\\=
\varepsilon^{-1}m_0(\Psi)\varphi(0)+m_1(\Psi)\varphi'(0)+O(\varepsilon)
\end{multline*}
as $\varepsilon \to 0$. Hence sequence $\langle\Psi_\varepsilon,\varphi\rangle$ converges for all $\varphi\in C_0^\infty(\mathbb{R})$ to a finite limit  if and only if $m_0(\Psi)=0$.
Moreover if $m_1(\Psi)\neq 0$, then the limit is nontrivial and the limit functional is defined by  $\Psi_0=-m_1(\Psi)\delta'(x)$, since $\langle\Psi_0,\varphi\rangle=m_1(\Psi)\varphi'(0)$.
\end{proof}

\begin{figure}[h]
\centering
  \includegraphics[scale=.6]{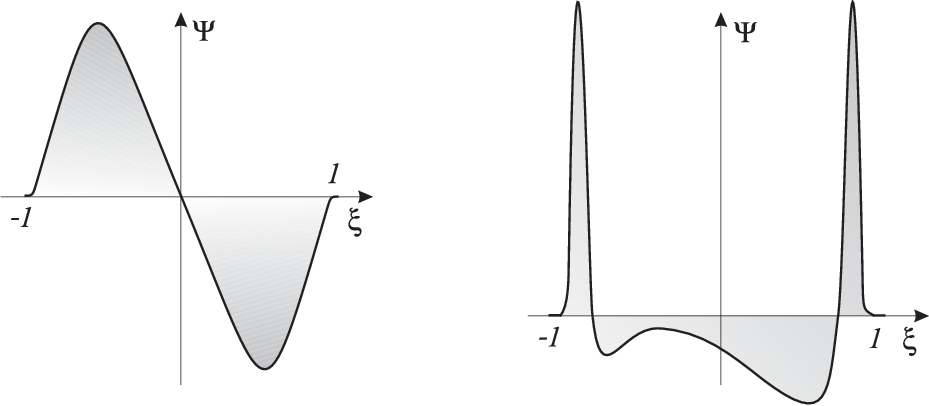}
  \caption{Plots of $\delta'$-like profiles}\label{FigPotential}
\end{figure}

The first condition in \eqref{CondDeltPrime} implies that function $\Psi$ must change a sign on $[-1,1]$. The second one involves an asymmetric property of  this function: if $\Psi(\xi)$ satisfies \eqref{CondDeltPrime}, then $\Psi(\xi-a)$ cannot be an even function for any $a$.
We introduce the set of \emph{$\delta'$-like potentials}
$$\mathcal{P}_0=\{\Psi\in \mathcal{P}\colon m_0(\Psi)=0,\, m_1(\Psi)=-1 \}.$$

Recall that $U(x)\to +\infty$ as $\abs{x}\to +\infty$. Then for all $\varepsilon>0$ the spectrum of $\mathcal{H}_{\varepsilon}(\alpha,\Psi)$ is real, discrete and simple.
Let $\{\lambda_k^\varepsilon(\alpha, \Psi)\}_{k=1}^\infty$ be the eigenva\-lues of  $\mathcal{H}_{\varepsilon}(\alpha,\Psi)$ enumerated in increasing. Suppose that $\{y_k^\varepsilon(x;\alpha, \Psi)\}_{k=1}^\infty$ is the orthonormal in $L_2(\mathbb{R})$ system of eigenfunctions.
\begin{g_theorem}\label{TheoremNeg}
For each pairs $(\alpha,\Psi)\in \mathbb{R}\times \mathcal{P}$  eigenvalues  $\lambda_k^\varepsilon(\alpha, \Psi)$ are continuous functions with respect to $\varepsilon \in (0,1)$.
Moreover, the eigenvalues are bounded from above as $\varepsilon\to 0$.
In the case of the $\delta'$-like potential, the spectrum of $\mathcal{H_{\varepsilon}}(\alpha,\Psi)$ is unbounded from below as $\varepsilon\to 0$: if $\alpha\neq 0$ and $\Psi\in \mathcal{P}_0$, then $\lambda_1^\varepsilon(\alpha, \Psi)\leq -\gamma\varepsilon^{-2}$
for some positive constant $\gamma$. There is only finite number $N(\alpha,\Psi)$ of eigenvalues that converge to $-\infty$ as $\varepsilon\to 0$.
\end{g_theorem}
\begin{proof} Fix $\alpha$ and $\Psi$. The quadratic form
$$
a_{\varepsilon}[u]=\int_{\mathbb{R}}\bigl(|u'|^2+(U+\alpha\Psi_\varepsilon)|u|^2\bigr)\,dx,\qquad u\in C_0^\infty(\mathbb{R})
$$
is a continuous function with respect to $\varepsilon\in(0,1)$.
The minimax principle
$$
\lambda_k^\varepsilon(\alpha, \Psi)=\inf\limits_{E_k}\,\sup\limits_{v\in
E_k,\: \norm{v}=1}a_{\varepsilon}[v]
$$
yields the continuity of eigenvalues  with respect to $\varepsilon$.
Here $E_k$ is a $k$-dimensional linear subspace of $C_0^\infty(\mathbb{R})$, and $\|\cdot\|$ denotes the norm in $L_2(\mathbb{R})$.

We choose subspace $E_k^*$, whose elements vanish in a neighborhood  of the origin.
Then, we have
\begin{equation}\label{InEqMinMax}
\lambda_k^\varepsilon(\alpha, \Psi)\leq\sup\limits_{v\in E_k^*,\:
\norm{v}=1}a_{\varepsilon}[v].
\end{equation}
For all sufficiently small $\varepsilon$ the restriction of $a_{\varepsilon}$ to the finite-dimensional space $E_k^*$ does not depend on $\varepsilon$. This yields the boundedness of the eigenvalues from above.

From now on,  $\Psi$ is a $\delta'$-like profile.
It is well known  \cite[p.~338]{ReedSimonIV} that the Schr\"{o}dinger operator $-\frac{d^2}{d\xi^2}+\alpha\Psi(\xi)$
on the line has  a negative eigenvalue  for all $\alpha\neq 0$ if and only if $m_0(\Psi)=0$.
Hence there exists  function $u$, $\|u\|=1$, such that
\begin{equation*}
    \mu=\int\limits_\mathbb R \bigl(u'^2+\alpha\Psi(\xi)u^2\bigr)\,d\xi<0.
\end{equation*}
Let $\zeta\in C^\infty_0(\mathbb R)$ be the positive  cutoff function such that  $\supp \zeta =[-2,2]$
and $\zeta(x)=1$ for $x\in [-1,1]$.
For abbreviation, we write $\zeta_\varepsilon$ and $\zeta'_\varepsilon$ instead of  $\zeta(\varepsilon\xi)$ and
$\zeta'(\varepsilon\xi)$ respectively. Let us consider  sequence
$u_\varepsilon(x)=\break\varepsilon^{-1/2}\zeta(x)u(\varepsilon^{-1}x)$. Then
 \begin{equation*}
\norm{u_\varepsilon}^2=\frac{1}{\varepsilon}\int\limits_\mathbb
R\zeta^2(x)u^2({\textstyle\frac{x}{\varepsilon}})\,dx= \int\limits_\mathbb
R\zeta^2_\varepsilon( \xi)u^2(\xi)\,d\xi\rightarrow \int\limits_\mathbb R u^2(\xi)\,d\xi=1
 \end{equation*}
as $\varepsilon\to 0$. Next, we conclude from the minimax principle  that
\begin{equation*}
    \lambda_1(\varepsilon,\alpha)=\inf\limits_{\substack{v\in
C_0^\infty,\\ \norm{v}=1}}a_{\varepsilon}[v]\leq
\frac{a_{\varepsilon}[u_\varepsilon]}{\|u_\varepsilon\|},
\end{equation*}
hence that
\begin{multline*}
\varepsilon^2\|u_\varepsilon\|\,\lambda_1(\varepsilon,\alpha)\leq
\varepsilon^2a_{\varepsilon}[u_\varepsilon]\\
= \varepsilon
\int\limits_\mathbb R\bigl((\zeta(x)u(\varepsilon^{-1}x))'^2+(U+\alpha \Psi_\varepsilon)(\zeta(x)u(\varepsilon^{-1}x))^2\bigr)\,dx\\
=\int\limits_\mathbb R\zeta^2_\varepsilon\bigl(u'^2+\alpha \Psi u^2\bigr)\,d\xi+
2\varepsilon\int\limits_\mathbb R\zeta_\varepsilon\zeta_\varepsilon' u u'\,d\xi\\+
\varepsilon^2\int\limits_\mathbb R\zeta_\varepsilon'^2u^2\,d\xi+
\varepsilon\int\limits_{-2}^2U(x)\zeta^2(x)u^2(\varepsilon^{-1}x)\,dx.
\end{multline*}
It is easily seen that the first integral converges to negative number  $\mu$ as $\varepsilon\to 0$, and  other terms go to zero. Thus for $\varepsilon$  sufficiently small the estimate
 $\varepsilon^2\|u_\varepsilon\|\,\lambda_1(\varepsilon,\alpha)\leq
\mu/2$ holds, and so  $\lambda_1(\varepsilon,\alpha)\leq -\gamma\varepsilon^{-2}$ for some  $\gamma>0$.

Let $N_\varepsilon^-(\alpha,\Psi)$ be the number of negative eigenvalues of $\mathcal{H}_{\varepsilon}(\alpha,\Psi)$.
Obviously $N(\alpha,\Psi)\leq N_\varepsilon^-(\alpha,\Psi)$ for $\varepsilon>0$ sufficiently small.
In the case of the continuous potential the inequality
\begin{equation}\label{Nminus}
   N_\varepsilon^-(\alpha,\Psi)\leq 1+\int\limits_{\mathbb R} |x|\,|U^-(x)|\,dx+|\alpha|\int\limits_{\mathbb R} |x|\,|\Psi^-_\varepsilon(x)|\,dx
\end{equation}
holds \cite[p.97]{BS}, where $f^-(x)=\min\{f(x), 0\}$ is the negative part of $f$.
Since $U$
increases as $x\to \infty$, $U^-$ is a function with compact support, and consequently
the first integral in \eqref{Nminus} is defined.
The support of $\Psi^-_\varepsilon$ lies in $[-\varepsilon,\varepsilon]$, we thus deduce
\begin{equation*}
\int\limits_{\mathbb R}
|x|\,|\Psi^-_\varepsilon(x)|\,dx=\varepsilon^{-2}\int\limits_{-\varepsilon}^{\varepsilon}
|x|\,|\Psi^-(\frac{x}{\varepsilon})|\,dx= \int\limits_{-1}^{1} |\xi|\,|\Psi^-(\xi)|\,d\xi.
\end{equation*}
In view of \eqref{Nminus} and the formula above
$N(\alpha,\Psi)$ is uniformly bounded with respect to  $\varepsilon$:
$$N(\alpha,\Psi)\leq N_\varepsilon^-(\alpha,\Psi)\leq c_1(U)+c_2(\Psi)|\alpha|$$
with constants $c_1(U)$, $c_2(\Psi)$  being positive.
\end{proof}

Therefore   the spectrum of $\mathcal{H_{\varepsilon}}(\alpha,\Psi)$ consists of two parts:
$\{\lambda_k^\varepsilon(\alpha, \Psi)\}_{k=1}^N$ is the set of eigenvalues that tend to $-\infty$ as $\varepsilon\to 0$, and
$\{\lambda_k^\varepsilon(\alpha, \Psi)\}_{k=N+1}^\infty$ is the set of bounded eigenvalues,
which we  call the bounded spectrum of $\mathcal{H_{\varepsilon}}(\alpha,\Psi)$. Notice that
eigenfunctions $y_1^\varepsilon(x;\alpha, \Psi), \dots, y_N^\varepsilon(x;\alpha, \Psi)$ have strongly oscillatory character as $\varepsilon\to 0$ in an neighborhood of the origin and vanish exponentially outside.
They converge to zero in $L_2(\mathbb{R})$ weakly.

\section{Asymptotics of Bounded Spectrum of  $\mathcal{H_{\varepsilon}}(\alpha,\Psi)$: Leading Terms}\label{SectionMainTerms}
Let us  consider the eigenvalue problem
\begin{equation}\label{MainProbl}
-y_\varepsilon''+\left(U(x)+\alpha \varepsilon^{-2}\Psi(\varepsilon^{-1}x)\right)y_\varepsilon=\lambda^\varepsilon y_\varepsilon,\qquad y_\varepsilon\in L_2(\mathbb{R}).
\end{equation}
Fix  eigenvalue $\lambda_k^\varepsilon(\alpha, \Psi)$ with  number $k>N(\alpha,\Psi)$. We denote it briefly by $\lambda^\varepsilon$.
Let $y_\varepsilon$ be the corresponding eigenfunction. We postulate asymptotic expansions for the eigenvalue
\begin{align}
\label{ExpanEVLow}
    \lambda^\varepsilon & \sim \lambda+\varepsilon\lambda_1+\varepsilon^2\lambda_2+\cdots,\\
\intertext{and two-scale expansions for the eigenfunction}
   \label{ExpanEFVLow}
  y_\varepsilon(x) & \sim v(x)+\varepsilon \,v_1(x)+\varepsilon^2 v_2(x)+\cdots\qquad
  \text{for}\; \abs{x}>\varepsilon,
 \\ y_\varepsilon(x) & \sim w(\varepsilon^{-1}x)+\varepsilon\, w_1(\varepsilon^{-1}x)+\varepsilon^2 w_2(\varepsilon^{-1}x)+
   \cdots\qquad   \text{for}\;\abs{x}\leq\varepsilon.\label{ExpanEFWLow}
\end{align}
Here  $v$, $v_i$ are defined for $x\in\mathbb{R}\setminus \{0\}$ and belong to $L_2(\mathbb{R})$.  Set $\xi=\varepsilon^{-1}x$. Functions $w$, $w_i$ are defined for $\xi\in [-1,1]$.
We also assume that $v$ is deferent from zero.
Series \eqref{ExpanEFVLow}, \eqref{ExpanEFWLow}  satisfy
the  coupling conditions
\begin{equation}\label{CondCoupl}
\left[y_\varepsilon\right]_{x=\pm\varepsilon}=0,\quad\bigl[y_\varepsilon'\bigr]_{x=\pm\varepsilon}=0
\end{equation}
at  points  $x=\pm \varepsilon$, where $[f]_{x=a}=f(a+0)-f(a-0)$ is the jump of $f$ at point $a$. This type of combined  asymptotic expansions has been used in particular in \cite{GolovatyCR}-\cite{GolovatySMZh}.

 Upon substituting  \eqref{ExpanEVLow}-\eqref{ExpanEFWLow} into equation \eqref{MainProbl}, we obtain
\begin{align}\label{VEquations}
  &- v''+U(x)v=\lambda v,\qquad x\in \mathbb{R}\setminus \{0\},\\ \label{V1Equations}
  &-v_1''+U(x)v_1=\lambda v_1+\lambda_1v, \qquad x\in \mathbb{R}\setminus \{0\},
\intertext{and also}
\label{WEquations}
&-w''+\alpha\Psi(\xi)w=0,\qquad \xi\in (-1,1), \\\label{W1Equations}
&-w_1''+\alpha\Psi(\xi)w_1=0,\qquad \xi\in (-1,1),\\\label{W2Equations}
&-w_2''+\alpha\Psi(\xi)w_2=\lambda w -U(0)w ,\qquad \xi\in(-1,1).
\end{align}
Next, substituting series \eqref{ExpanEFVLow}, \eqref{ExpanEFWLow} into the coupling conditions yields
\begin{align*}
        v(\pm\varepsilon)+\varepsilon \,v_1(\pm\varepsilon)+\cdots &\sim w(\pm1)+\varepsilon\, w_1(\pm1)+\cdots ,\\
     v'(\pm\varepsilon)+\varepsilon \,v_1'(\pm\varepsilon)+\cdots &\sim \varepsilon^{-1}w'(\pm1)+ w_1'(\pm1)+\varepsilon w_2'(\pm1)+\cdots.
\end{align*}
We can now expand $v(\pm\varepsilon)$, $v'(\pm\varepsilon)$ into the formal Taylor series
about $x = \pm0$. Then
\begin{gather}\label{AdCondVW}
v(-0)=w(-1),\qquad v(+0)=w(+1),\\
\label{AdCondWPrime}
w'(-1)=0, \qquad w'(1)=0, \\
\label{AdCondW1Prime}
v'(-0)=w_1'(-1),\qquad v'(+0)=w_1'(1),\\
\label{AdCondVW1}
v_1(-0)- v '(-0)=w_1(-1), \quad v_1(+0)+ v '(+0)=w_1(1),\\\label{AdCondV1W2}
v_1'(-0)-v ''(-0)=w_2'(-1),\quad v_1'(+0)+v ''(+0)=w_2'(1).
\end{gather}

 It follows that $v$  satisfies equation \eqref{VEquations} on each half-axis and $w$ is a solution to the problem
\begin{equation}\label{W0Problem}
-w''+\alpha \Psi w=0,\quad\xi\in(-1,1),\qquad
w'(-1)=0,\; w'(1)=0.
\end{equation}
Moreover  both these functions  satisfy  coupling conditions \eqref{AdCondVW}.
Problem \eqref{W0Problem} is decisive in our next considerations, because it
contains information about a kind of the singular perturbation, namely  the profile $\Psi$ and the coupling constant $\alpha$.
The first and primary question is whether there  exists its nontrivial solution.

\subsection{Resonant Set and Coupling Function}
Problem \eqref{W0Problem} can be viewed as a spectral problem with  spectral parameter $\alpha$.
Since function $\Psi$ is sign-changing, it is natural to introduce a space with an indefinite metric.
Let $\mathcal{K}$ be the weight $L_2$-space with the scalar product
$(f,g)=\int\limits_{-1}^1\abs{\Psi}f\overline{g}\,d\xi$.
Let us introduce in $\mathcal{K}$ the indefinite metric $[f,g]=\int\limits_{-1}^1\Psi f \overline{g}\,d\xi$.
Then the pair $(\mathcal{K},[\cdot,\cdot])$  is called a \emph{Krein space}
(see \cite{IA} for the original definition). In the Krein space there exists the \emph{fundamental symmetry}
$Jf=\sgn\Psi\, \cdot f$ such that $[f,g]=(Jf,g)$ for all $f,g\in \mathcal{K}$.

Let $T$ be a closed densely defined operator in $\mathcal{K}$.
The $J$-\textit{adjoint} operator of $T$ is defined by the
relation $[Tx,y]=[x,T^{[*]}y]$ for $x\in \mathcal{D}(T)$
on the set of all $y\in \mathcal{K}$ such that the mapping $x\mapsto [Tx, y]$ is a continuous linear functional on
$\mathcal{D}(T)$. The operator $T$ is called $J$\textit{-selfadjoint} if $T = T^{[*]}$.
The operator  $T$ is called \emph{$J$-nonnegative} if $[Tx,x]\geq0$ for  $x\in \mathcal{D}(T)$.
\begin{g_lemma}\label{LemmaJ}
Let $T$ be a $J$-selfadjoint and $J$-nonnegative operator with a nonempty resolvent set.
Spectrum $\sigma(T)$ is real, moreover, the residual spectrum is empty. If $\lambda\neq 0$ is
an eigenvalue of $T$, then it is simple (e.i., $\ker (T-\lambda)=\ker (T-\lambda)^2$).
If $0$ is an eigenvalue of $T$,  then its Riesz index is equal or less than $2$, e.i., $\ker T^2=\ker T^3$.
Generally, $0$ is not a simple eigenvalue \cite[p. 138]{IA}.
\end{g_lemma}
We introduce in $\mathcal{K}$ the operator $\mathcal{T}_\Psi=-\frac{1}{\Psi(\xi)}\frac{d^2}{d\xi^2}$ with the domain
$$
\mathcal{D}(\mathcal{T}_\Psi)=\{\bigl.f\in \mathcal{K}\:\bigr|\: f\in
W_2^2(-1,1),\; \Psi^{-1}f''\in \mathcal{K},\; f'(-1)=0,\:f'(1)=0
\}.
$$
Spectral equation $\mathcal{T}_\Psi w=-\alpha w$ is associated with \eqref{W0Problem}.
Differential equations with the indefinite weights have been widely studied  \cite{CurJDE, IohvidovKreinLanger, KarabashKostenkoMalamud, PiatkovSMZh}. For the sake of completeness we describe the main properties of $\mathcal{T}_\Psi$.
\begin{g_theorem}\label{ThJproperties}
For each  $\Psi\in \mathcal{P}$  operator $\mathcal{T}_\Psi$ is $J$-selfadjoint and $J$-non\-ne\-ga\-tive.
\end{g_theorem}
\begin{proof}
Given $f\in \mathcal{D}(\mathcal{T}_\Psi)$, we have
\begin{equation*}
[\mathcal{T}_\Psi f,g]=-\int\limits_{-1}^1f''\overline{g}\,d\xi=
f(1)\,\overline{g'(1)}-f(-1)\,\overline{g'(-1)}
-\int\limits_{-1}^1f\,\overline{{g}''}\,d\xi.
\end{equation*}
Then  identity $[\mathcal{T}_\Psi f,g]=[f,\mathcal{T}_\Psi^{[*]} g]$ holds
if $g\in \mathcal{D}(\mathcal{T}_\Psi)$ and $\mathcal{T}_\Psi^{[*]} g=-\Psi^{-1}g''$. Hence  $\mathcal{T}_\Psi$ is $J$-selfadjoint. Next, for all $f\in \mathcal{D}(\mathcal{T}_\Psi)$
$$
[\mathcal{T}_\Psi f,f]=-\int\limits_{-1}^1f'' \overline{f}\,d\xi=
\int\limits_{-1}^1\abs{f'}^2\,d\xi\geq0,
$$
so that   $\mathcal{T}_\Psi$ is a $J$-nonnegative operator.
\end{proof}

\begin{g_theorem}\label{ThSpectrumOfT}%
\begin{itemize}%\itemindent20pt
\item[$(i)$] Given $\Psi\in \mathcal{P}$, the spectrum of  $\mathcal{T}_\Psi$ is real and discrete.
All nonzero eigenvalues are simple. If $m_0(\Psi)\neq0$, then $0$ is a simple eigenvalue.
In the general case, $\ker \mathcal{T}_\Psi\neq\ker \mathcal{T}_\Psi^2$.

\item[$(ii)$] If $\Psi\in \mathcal{P}_0$, then the spectrum of $\mathcal{T}_\Psi$ has two accumulation points $-\infty$ and $+\infty$. Moreover, $\ker \mathcal{T}_\Psi\neq \ker \mathcal{T}_\Psi^2$.
\end{itemize}
\end{g_theorem}
\begin{proof}
To start with, show that the resolvent set of $\mathcal{T}_\Psi$ is nonempty.
Each solution of the homogenous problem
\begin{equation}
\label{HomogProb}
g''+i\Psi g=0,\quad\xi\in(-1,1),\qquad
g'(-1)=0,\; g'(1)=0,
\end{equation}
satisfies equality
$$
\int\limits_{-1}^1\abs{g'}^2d\xi-i\int\limits_{-1}^1\Psi\abs{g}^2d\xi=0.
$$
Since function $\Psi$ is real-valued, we deduce that $g$ is a constant function.
Obviously, the zero function can only be a solution of (\ref{HomogProb}).
Then the nonhomogeneous problem $g''+i\Psi g=f$, $g'(-1)=0$, $ g'(1)=0$ admits a unique solution
for each $f\in \mathcal{K}$. Hence $i\in \rho(\mathcal{T}_\Psi)$ and $g=R(i,\mathcal{T}_\Psi)f$, where $R(\lambda,\mathcal{T}_\Psi)$ is the resolvent of  $\mathcal{T}_\Psi$.

Therefore spectrum  $\sigma(\mathcal{T}_\Psi)$ is real by
Lemma~\ref{LemmaJ} and Theorem~\ref{ThJproperties}. The compactness of $R(i,\mathcal{T}_\Psi)$ follows from
the sequence of embeddings
$$
\mathcal{D}(\mathcal{T}_\Psi)\subset W_2^2(-1,1)\hookrightarrow L_2(-1,1)\subset
\mathcal{K}.
$$
As a consequence we have $\sigma(\mathcal{T}_\Psi)=\sigma_p(\mathcal{T}_\Psi)$.

Each eigenvalue of $\mathcal{T}_\Psi$ is simple. Indeed, if there exist two linear independent
eigenfunctions $\varphi$ and $\psi$ corresponding to eigenvalue $\alpha$, then function
$w(\xi)=\psi(1)\varphi(\xi)-\varphi(1)\psi(\xi)$ is nonzero since both values
$\varphi(1)$, $\psi(1)$ are different from $0$.
But this is impossible, because $w=0$ as a unique solution to the Cauchy problem $w''+\alpha\Psi w=0$, $w(1)=0$, $w'(1)=0$.

Assume now $\Psi\in \mathcal{P}_0$. Since potential $\Psi$ changes a sign,  the spectrum  $\sigma(\mathcal{T}_\Psi)$ is unbounded  in both directions \cite{CurJDE} (see Remark~\ref{RemarkOnOddPsi}).
Moreover  invariant subspace  $\ker \mathcal{T}_\Psi^2$ is generated by eigenvector $w=1$ and root vector
$$
w_*(\xi)=\int\limits_{-1}^{\xi}(t-\xi)\Psi(t)\,dt,
$$
which is a solution to the problem $w_*''=-\Psi(\xi)$, $w'_*(-1)=0$, $w'_*(1)=0$.
This problem admits a solution if and only if $m_0(\Psi)=0$. In view of Lemma~\ref{LemmaJ} there are no linear independent root  vectors other than $w_*$.
\end{proof}

\begin{g_remark}\label{RemarkOnOddPsi}\rm
The spectrum of $\mathcal{T}_\Psi$ is a symmetric set  with respect to the origin for each \emph{odd} profile $\Psi$:
if $\alpha$ is an eigenvalue with eigenfunction $w(\xi)$, then $-\alpha$ is also an eigenvalue with eigenfunction $w(-\xi)$.
\end{g_remark}

We introduce the set $\Sigma_\Psi=\{\alpha\in \mathbb{R}\colon -\alpha\in \sigma(\mathcal{T}_\Psi)\} $. We call
$\Sigma_\Psi$ the \emph{resonant set} of  profile  $\Psi$.
Problem \eqref{W0Problem} has a nontrivial solutions if and only if the coupling constant $\alpha$ belongs to  $\Sigma_\Psi$. Suppose $W_\alpha$ is an eigenfunction of \eqref{W0Problem} with eigenvalue $\alpha\in \Sigma_\Psi$. Thus values $W_\alpha(-1)$ and $W_\alpha(1)$ are nonzero and so the quotient
$$
\theta_\Psi (\alpha)=\frac{W_\alpha(1)}{W_\alpha(-1)}
$$
is correctly defined.
Quantity $\theta_\Psi (\alpha)$ does not depend on the choice of an eigenfunction and can be considered as  function $\theta_\Psi \colon\Sigma_\Psi \to \mathbb{R}$ defined at points of the resonant set. We call $\theta_\Psi$
the  \emph{coupling function} of  $\Psi$. Hence each profile $\Psi$ generates the resonant set and the coupling function.

\subsection{Limit Operator $\mathcal{H}(\alpha,\Psi)$}\label{SubsecLimitOperator}
 We now continue to construct the asymptotics. We distinguish two different cases.
First assume  \emph{$\alpha$ does not belong to  resonant set $\Sigma_\Psi$.}
Then problem \eqref{W0Problem} has  a trivial solution $w=0$ only, and so from  coupling conditions \eqref{AdCondVW} we deduce $v(-0)=v(+0)=0$.  Recalling \eqref{VEquations}, we obtain
\begin{equation}\label{LimitProblemW=0}
 - v''+Uv=\lambda v,\quad x\in\mathbb{R}\setminus\{0\}, \qquad v(0)=0, \quad v\in L_2(\mathbb{R}).
\end{equation}
Function $v$  is nonzero by our assumption, thus $\lambda$ is an eigenvalue of the problem \eqref{LimitProblemW=0}.
Therefore we can introduce limit operator $\mathcal{H}(\alpha,\Psi)$, which is the direct sum $S_-\oplus S_+$
of the Schr\"{o}dinger operators acting on half-axes. Operators $S_-$ and $S_+$ are associated with problems
\begin{equation}\label{AxisProblem}
\begin{cases}
   -v''+Uv=\lambda v,\;\; x\in \mathbb{R}_-,\\
   \phantom{-}v(0)=0, \quad v\in L_2(\mathbb{R}_-),
\end{cases}\qquad
\begin{cases}
   -v''+Uv=\lambda v,\;\; x\in \mathbb{R}_+,\\
   \:v(0)=0, \quad v\in L_2(\mathbb{R}_+)
\end{cases}
\end{equation}
respectively, where $\mathbb{R}_-$ and $\mathbb{R}_+$ are the negative and positive half-axes.
Spectra of  $S_-$ and $S_+$ are real, discrete and simple. Obviously, $$\sigma(\mathcal{H}(\alpha,\Psi))=\sigma(S_-)\cup\sigma(S_+).$$

Let us now suppose that $\alpha$ \emph{belongs to  resonant set} $\Sigma_\Psi$. We set $w=aW_\alpha$, where $W_\alpha$ is the corresponding eigenfunction of \eqref{W0Problem}, and $a\neq 0$.
We conclude from \eqref{AdCondVW} that  $v(-0)=aW_\alpha(-1)$, $v(+0)=aW_\alpha(1)$, hence that
\begin{equation}\label{ThetaVCond}
    v(+0)=\theta_\Psi (\alpha)v(-0).
\end{equation}
According to \eqref{W1Equations}, \eqref{AdCondW1Prime}  the next term $w_1$ of series \eqref{ExpanEFWLow} can be found by solving the problem
\begin{equation}
\label{W1Problem}
-w_1''+\alpha\Psi w_1=0,\:\xi\in(-1,1),\quad
w_1'(-1)=v'(-0),\: w_1'(1)=v'(+0).
\end{equation}
Because $\alpha$ is an eigenvalue of \eqref{W0Problem}, the problem admits a solution if and only if $W_\alpha(1)v'(+0)=W_\alpha(-1)v'(-0)$.
This solvability condition can be derived multiplying the equation by $W_\alpha$ and integrating
twice by parts.
It may be written in the form
\begin{equation}\label{ThetaVPrimeCond}
    \theta_\Psi (\alpha)v'(+0)=v'(-0).
\end{equation}
Combining  \eqref{VEquations}, \eqref{ThetaVCond} and \eqref{ThetaVPrimeCond} we deduce that $v$ must be an eigenfunction of the  problem
\begin{equation}\label{LimitProblemV}
\begin{cases}
    -v''+U v=\lambda v,\qquad x\in\mathbb{R}\setminus\{0\},\\
    v(+0)-\theta_\Psi (\alpha)v(-0)=0,\quad \theta_\Psi (\alpha) v'(+0)-v'(-0)=0.
\end{cases}
  \end{equation}
Thus the connected selfadjoint extension of operator $L$
\begin{equation}\label{LimitOperatorResonance}
\begin{gathered}
\hskip-6.5cm    \mathcal{H}(\alpha,\Psi)=-\frac{d^2}{dx^2}+U(x),\\
\mathcal{D}(\mathcal{H}(\alpha,\Psi)) =\{f\in \mathcal{D}(L^*)\colon
f(+0)\hskip5cm\\\hskip3.5cm=\theta_\Psi (\alpha)f(-0),\: \theta_\Psi (\alpha)f'(+0)=f'(-0)
\}
\end{gathered}
\end{equation}
is associated with \eqref{LimitProblemV}.

\begin{g_remark}\rm \label{RemarkV(0)=0} In the case $\alpha\in \Sigma_\Psi$ we  can also assume that
 $a=0$, that is, $w=0$.  Then $v(0)=0$, and so  $v$ is an eigenfunction of $S_-\oplus S_+$. But condition \eqref{ThetaVPrimeCond}
has  nevertheless to hold true; for otherwise we could not solve \eqref{W1Problem}.
If $v$ is an eigenfunction corresponding to a simple eigenvalue of $S_-\oplus S_+$, then at least one of values $v'(-0)$ and $v'(+0)$ equals $0$, since $v$ vanishes on at least one of half-axes. Then condition  \eqref{ThetaVPrimeCond}
yields $v'(-0)=v'(+0)=0$, which is impossible.
For a double eigenvalue of  $S_-\oplus S_+$, there exists a unique line in the 2-dimensional eigenspace, for which condition \eqref{ThetaVPrimeCond} can be satisfied.
Vectors of this line are also eigenfunctions of \eqref{LimitProblemV}, since condition \eqref{ThetaVCond}
is trivially valid. Therefore in the resonant case  the leading term  $w$ of \eqref{ExpanEFWLow} is zero if and only if  $\lambda\in\sigma(\mathcal{H}(\alpha,\Psi))\cap\sigma(S_-\oplus S_+)$.
This intersection consists of all double eigenvalues of $S_-\oplus S_+$. Note that for operator $\mathcal{H}(\alpha,\Psi)$ these eigenvalues are simple.
\end{g_remark}

\begin{g_remark}\rm For $\alpha=0$ operators $\mathcal{H_{\varepsilon}}(0,\Psi)$ don't depend on $\varepsilon$ and coincide with the unperturbed  Schr\"{o}dinger operator with potential $U$. On the other hand, the resonant set always contains
 $0$, and so operator $\mathcal{H}(0,\Psi)$ must be defined as in \eqref{LimitOperatorResonance}.
Obviously $\theta_\Psi (0)=1$, because the constant eigenfunction corresponds to the zero eigenvalue of
$\mathcal{T}_\Psi$. It follows that  $\mathcal{H_{\varepsilon}}(0,\Psi)=\mathcal{H}(0,\Psi)$ for all $\varepsilon>0$.
\end{g_remark}

After all, for each profile $\Psi$ we have defined the family of selfadjoint operators  $\{\mathcal{H}(\alpha,\Psi)\}_{\alpha\in\mathbb{R}}$.
If coupling constant $\alpha$ doesn't belong to resonant set $\Sigma_\Psi$, then operator
$\mathcal{H}(\alpha,\Psi)$ is a direct sum of the Schr\"{o}dinger operators acting on half-axes. The corresponding eigenfunctions describe  states of a quantum system, when a particle can be found with probability $1$ on one of the half-axes. In this case we obtain what will be referred to as the \emph{closed $\delta'$-barrier}.
 If $\alpha\in \Sigma_\Psi$, then operator $\mathcal{H}(\alpha,\Psi)$ given by \eqref{LimitOperatorResonance} is the  connected  selfadjoint extension with  matrix
\begin{equation}\label{MatrixC}
C_\alpha=\begin{pmatrix} \theta_\Psi (\alpha) & 0 \\ 0 & \theta_\Psi (\alpha)^{-1} \end{pmatrix}
\end{equation}
in \eqref{ConnectedCond}. For resonant coupling constants  a particle can permeate through the barrier with some probability  (the case of an \emph{open $\delta'$-barrier}).

\section{What Is the $\delta'$-Potential?}\label{SectionHistory}
\subsection{Historical Remarks}
The formal Hamiltonians with the Dirac functions and their derivatives in potentials were studied since the eighties of last century. We  refer, for example, to  research monograph \cite{Albeverio2edition} (the first edition, 1988) and papers \cite{SebaCzechJPhys86, SebRMP, Grossmann}.
Incidentally, in solvable models one distinguishes  two different phenomena dealing with the $\delta'$-function:
\textit{the $\delta'$-interaction} and  \textit{the point dipole interaction (the $\delta'$-potential).}

The first result about the $\delta'$-interaction is due to S.~Albeverio, F.~Geszte\-sy, R.~H{\o}egh-Krohn and  H.~Holden \cite{Albeverio2edition}.
The standard definition  is the one-parameter family of selfadjoint extensions
\begin{equation}\label{AlbeverioDefn}
\begin{gathered}
A_\beta=-\frac{d^2}{dx^2},\\ \mathcal{D}(A_\beta)=\big\{f\in
W_2^2(\mathbb{R}\setminus\{0\})\colon f'(-0)=f'(+0),\hskip3cm\\\hskip4.5cm
f(+0)-f(-0)=\beta f'(0)\big\}.
\end{gathered}
\end{equation}
P.~\v{S}eba \cite{SebRMP} showed   that these selfadjoint extensions correspond to the heuristic operator
\begin{equation*}
    -\frac{d^2}{dx^2}+\beta|\delta'(x)\,\rangle\,\langle\,\delta'(x)|,
\end{equation*}
where the rank one perturbation $|f\,\rangle\,\langle\,f|$ is defined by
$$
(|f\,\rangle\,\langle\,f|\varphi\,\rangle)(x)=f(x)\int\limits_\mathbb R
f(y)\varphi(y)\,dy.
$$
In the same paper  P.~\v{S}eba  also studied  formal Hamiltonian
$-\frac{d^2}{dx^2}+\alpha\delta'(x)$, where  $\delta'(x)$ was considered as the limit in the topology of  $\mathcal{D}'(\mathbb R)$ of  linear combination
$$
\frac{1}{2\varepsilon}\bigr(\delta(x+\varepsilon)-\delta(x-\varepsilon)\bigl).
$$
The limit operator is referred to as  the $\delta$-interaction dipole Hamiltonian.
The author proved that this Hamiltonian is the direct sum $S_-\oplus S_+$ with $U=0$ (see \eqref{AxisProblem}).
For instance, in Theorem~4 one proves that the family of the  Schr\"{o}dinger operators
with smooth potentials
\begin{equation}\label{SebaOperators}
    S_\varepsilon(V)= -\frac{d^2}{dx^2}+\frac{1}
    {\varepsilon^2}V\Bigl(\frac{x}{\varepsilon}\Bigr),\qquad x\in\mathbb R
\end{equation}
converges to $S_-\oplus S_+$ in the norm resolvent sense, provided potential $V\in C^\infty_0(\mathbb
R)$ fulfils condition  $m_0(V)=0$. From viewpoint of the scattering theory this means that the  $\delta'$ is a totally reflecting wall.

But these conclusions are at variance with the results presented in \cite{ChristianZolotarIermak03,ToyamaNogami,
Zolotaryuk06, Zolotaryuk08},  where the scattering properties of  $\delta'$-like potentials constructed
in the form of squeezed rectangles were studied in the zero-range limit
(a simple such problem is treated in  Subsection~\ref{subsecScattering}).
The resonances in the transmission probability were firstly described  in \cite{ChristianZolotarIermak03}.
The authors  obtained a discrete set of resonant coupling constants $\alpha_n$
for which there is a partial transmission decreasing as $\alpha_n$ is larger. These values are roots of a
transcendent equation, which depends on the regularization, that is to say, the profile of  piecewise constant potentials.

It is worth to  point out preprint \cite{KurElaMSI93}.
P.~Kurasov and N.~Elander proposed to define the Hamiltonian with the $\delta'$-potential as the second
derivative operator in $L_2(\mathbb R)$ subject to the coupling conditions
\begin{equation}\label{KurasovDefn}
\begin{gathered}
f(+0)-f(-0)=\frac{\alpha}{2}(f(+0)+f(-0)), \\[5pt]
f'(+0)-f'(-0)=-\frac{\alpha}{2}(f'(+0)+f'(-0)).
\end{gathered}
\end{equation}
This definition is based on the extension of Dirac's functions
to the case of discontinuous test functions \cite{KurasovJMAA96}:
$$\langle\delta^{(n)}(x),\varphi(x)
\rangle=\frac{(-1)^{n}}{2}(\varphi^{(n)}(+0)+\varphi^{(n)}(-0)).$$
The same definition was introduced by L.~Nizhnik \cite{NizhFAA2006}, who studied the Schr\"{o}dinger operator with
the $\delta'$-potential in Sobolev space $W_2^3(\mathbb R\setminus 0)$.

Alternative approaches to the problem can be found in \cite{GriffithsJPhA1993, GesHolJPA, KiselevJMAA1997}.
In \cite[p.~339]{AlbeverioKurasov} we read that the relation between the set of point interactions (selfadjoint extensions from Lemma~\ref{AllExtensionsLemma}) and the set of interactions defined by the $\delta'$-potential can not be established  without additional assumptions on the symmetry properties of the interaction.
Different assumptions can lead to different one dimensional families of point interactions corresponding to the formal expression $-\frac{d^2}{dx^2}+\alpha\delta'(x)$.
It favours the view that this problem contains  hidden parameters.

\subsection{Solvable Model for Hamiltonian  with Potential
$\alpha \varepsilon^{-2}\Psi(\varepsilon^{-1}x)$}\label{SubsecWhatIsDeltaPrime}

The formal asymptotic results obtained in the previous section point to the following fact:
the solvable model that is appropriate to operator $-\frac{d^2}{dx^2}+\alpha\delta'(x)$
is not specific for the given point interaction, it depends on the way in which  the derivative of Dirac's function is approximated in the weak topology. In fact,  profile $\Psi$ is a hidden parameter in the conventional formulation of problem how to define the $\delta'$-potential.

\emph{However, the solvable model that most closely corresponds to the motion of a quantum particle in potential $\frac{\alpha}{\varepsilon^{2}}\Psi(\varepsilon^{-1}x)$ is uniquely defined. The model is described by the family of operators $A(\alpha,\Psi)=-\frac{d^2}{dx^2}$ with domains
\begin{align*}
    &\mathcal{D}(A(\alpha,\Psi))=\{f\in W_2^2(\mathbb R\setminus 0)\colon f(-0)=f(+0)=0\} \quad\text{for } \alpha\not\in \Sigma_\Psi,\\
    &\mathcal{D}(A(\alpha,\Psi))=\{f\in W_2^2(\mathbb R\setminus 0)\colon f(+0)
=\theta_\Psi (\alpha)f(-0),\\ &\hskip5.2cm \theta_\Psi (\alpha)f'(+0)=f'(-0)\} \;\text{for
} \alpha\in \Sigma_\Psi,
\end{align*}
where $\Sigma_\Psi$ and $\theta_\Psi$ are the resonance set and the coupling function of
$\Psi$.}

This solvable model is in good agreement with the results obtained previously in
\cite{ChristianZolotarIermak03,ToyamaNogami, Zolotaryuk06, Zolotaryuk08}.
We actually describe the phenomenon of resonance in the transmission probability for arbitrary profiles $\Psi$
and obtain the resonance set as a spectral characteristic of the profile.

It is worth to compare our definition with that proposed in \cite{KurElaMSI93, NizhFAA2006} (see  \eqref{KurasovDefn}), in which the coupling matrix
\begin{equation*}\label{MatrixKurasov}
C_\alpha=\begin{pmatrix}\frac{2+\alpha}{2-\alpha} & 0\\ 0 & \frac{2-\alpha}{2+\alpha} \end{pmatrix}
\end{equation*}
is also diagonal provided $|\alpha|\neq 2$. In addition, for two values of $\alpha$ the model is described by
separated selfadjoint extensions, namely   $f'(-0)=0$,  $f(+0)=0$ for $\alpha=-2$, and
$f(-0)=0$, $f'(+0)=0$ for $\alpha=2$.

From the physical viewpoint  smooth potential $U$ has no effect on the coupling conditions at the origin.
Moreover,  the resonance set and the coupling function don't depend upon $U$.
We cannot but assume that $U=0$ in our computations above, because operators $\mathcal{H}_{\varepsilon}(\alpha,\Psi)$
lose discreteness of  the spectrum.  We shall give two examples of exactly solvable models with  piecewise constant $\delta'$-like potentials in order to motivate our definition of the point dipole interaction.

\subsection{Transition of Quantum Particle\\ Through $\delta'$-Like Potential}\label{subsecScattering}
Let us consider the equation
\begin{equation}\label{SolvModel2}
    -y''+\alpha\varepsilon^{-2}\Psi(\varepsilon^{-1}x)y=k^2 y,\quad x\in\mathbb{R},
\end{equation}
with $\Psi$ being a $\delta'$-like profile, such that  $\Psi(\xi)=1$ for $\xi\in(-1,0)$ and
$\Psi(\xi)=-1$ for $\xi\in(0,1)$, $k>0$. Clearly the results obtained above are still true for the piecewise constant potentials. Assume that $\alpha=\varkappa^2>0$. It is desired to find solution $y_\varepsilon(x;\varkappa,k)$ of \eqref{SolvModel2}, such that $y_\varepsilon(x;\varkappa,k)=e^{ikx}+R_\varepsilon(\varkappa,k)e^{-ikx}$ for
$x<-\varepsilon$ and $y_\varepsilon(x;\varkappa,k)=T_\varepsilon(\varkappa,k)e^{ikx}$ for
$x>\varepsilon$.
Quantities $|R_\varepsilon(\kappa,k)|^2$, $|T_\varepsilon(\kappa,k)|^2$ are called a \emph{reflection coefficient} and  a \emph{transmission coefficient} respectively. They have a probabilistic meaning in view of
$|R_\varepsilon(\varkappa,k)|^2+|T_\varepsilon(\varkappa,k)|^2=1$. The same example is treated in   \cite{ChristianZolotarIermak03}. We adduce this example to show the direct relationship between  scattering data and the coupling function introduced above.

For the profile  under consideration  resonant set  $\Sigma_\Psi$ is symmetric with respect to the origin, and
its positive part consists of the  roots of transcendent equation
$h(\sqrt{\alpha})=0$, where $h(\varkappa)=\varkappa(\tnh\varkappa-\tng\varkappa)$. The coupling function is given by
\begin{equation}\label{ThetaForSolvModel}
    \theta_\Psi (\alpha)=\frac{\csh\sqrt{\alpha}}{\cos\sqrt{\alpha}}\quad\text{if }\alpha\geq 0,\qquad
    \theta_\Psi (\alpha)=\frac{\cos\sqrt{-\alpha}}{\csh\sqrt{-\alpha}}\quad\text{if }\alpha<0.
\end{equation}

Let $u_1(\xi;\varkappa,\tau)$, $u_2(\xi;\varkappa,\tau)$ be the fundamental system of solutions of
$$
-u''+(\varkappa^2\Psi(\xi)-\tau^2)u= 0, \quad \xi\in(-1,1),
$$
such that its elements are smooth functions of $\varkappa$ and $\tau$.
Such solutions can be constructed explicitly with trigonometric and  hyperbolic functions.
Then we can find  solution  $y_\varepsilon(x;\varkappa,k)$ in the form
\begin{equation*}
y_\varepsilon(x;\varkappa,k)=
\begin{cases}
e^{ikx}+R_\varepsilon\,e^{-ikx}& \text{for } x<-\varepsilon,\\
C_{\varepsilon,1}\,u_1(\frac{x}{\varepsilon},\varkappa, \varepsilon k)+C_{\varepsilon,2}\,u_2(\frac{x}{\varepsilon},\varkappa, \varepsilon k) & \text{for } |x|<\varepsilon,\\
T_\varepsilon\,e^{ikx} & \text{for } x>\varepsilon.
\end{cases}
\end{equation*}
where we intend to select the coefficients to satisfy  the $C^1$-smoothness conditions at points $x=\pm \varepsilon$.
A direct computation gives
\begin{equation}\label{Teps}
  T_\varepsilon(\varkappa,k)= \frac{2i\varepsilon k e^{-2i\varepsilon k}}{(2i\varepsilon k-h(\varkappa))\cos\varkappa \csh\varkappa+O(\varepsilon^2k^2)}
\end{equation}
as $\varepsilon k\to 0$. The asymptotic behavior of the transmission coefficient as $\varepsilon\to 0$ depends upon whether $\alpha=\kappa^2$ belongs to the resonant set. Namely
\begin{align}\label{TepsAsymptotics}
  &|T_\varepsilon(\varkappa,k)|^2= \frac{4\varepsilon^2 k^2}{h^2(\varkappa)\cos^2\varkappa \csh^2\varkappa}\cdot
  \bigl(1+O\left(\varepsilon k\right)\bigr)\quad\text{for }
  \varkappa^2\not\in \Sigma_\Psi,\\\label{TepsAsymptoticsR}
  &|T_\varepsilon(\varkappa,k)|^2= \frac{1}{\cos^2\varkappa \csh^2\varkappa}\cdot\bigl(1+O(\varepsilon k)\bigr)\quad\text{for }
  \varkappa^2\in \Sigma_\Psi.
\end{align}
Whence it follows that the transmission coefficient goes to zero outside  the resonant set, and
it tends to a positive limit for $\alpha\in \Sigma_{\Psi}$. This limit $|T(\alpha)|^2$ is independent of $k$
and can be represented  in terms of the coupling function
\begin{equation*}
  |T(\alpha)|^2=\frac{4\theta_\Psi^2(\alpha)}{(1+\theta_\Psi^2(\alpha))^2}, \qquad \alpha\in \Sigma_\Psi.
\end{equation*}
This formula is still true  for negative $\alpha$. According to \eqref{ThetaForSolvModel},  transmission probability $|T(\alpha)|^2$ tends to zero as $|\alpha|\to +\infty$.

Therefore  operators $A(\alpha,\Psi)$, defined in the previous subsection, are  the \emph{connected} selfadjoint extensions if and only  if $\alpha\in \Sigma_\Psi$.

\subsection{Sturm-Liouville Operator With  $\delta'$-Like Potential}
Let $(a,b)$ be a bounded interval, which contains the origin. We consider the eigenvalue problem
\begin{equation}\label{SolvModel1}
    -y''+\alpha\varepsilon^{-2}\Psi(\varepsilon^{-1}x)y=\lambda y,\quad x\in(a,b),\qquad y(a)=0,\quad y(b)=0,
\end{equation}
where $\Psi$ is the piecewise constant $\delta'$-like potential as in the previous subsection.
The fundamental system of solutions of \eqref{SolvModel1}  can be constructed explicitly with trigonometric and  hyperbolic functions.

We  perform  the calculations for positive $\alpha$ only. Let us introduce the notation
$\omega=\sqrt{\lambda}$, $\varkappa=\sqrt{\alpha}$,
$h_1(\varkappa)=\tnh\varkappa\tng\varkappa-1$ and
$g(\varkappa)=(1+\tnh\varkappa\tng\varkappa)(1-\tnh\varkappa\tng\varkappa)^{-1}$.
The characteristic determinant $\Delta(\varepsilon, \varkappa;\omega)$ of \eqref{SolvModel1}, whose roots are eigenfrequencies $\omega_\varepsilon$, admits the asymptotic representation as $\varepsilon \omega\to 0$:
\begin{multline}\label{DetTeylor1}
\Delta(\varepsilon, \varkappa;\omega)= h(\varkappa)\bigl\{\tng a\omega \tng b\omega+\varepsilon\omega \,\bigl(\tng b\omega -\tng a\omega\bigr)\bigr\}\\
+\varepsilon\omega\, h_1(\varkappa)\bigl(\tng b\omega-g(\varkappa)\tng a\omega\bigr)
+O(\varepsilon^2\omega^2)
\end{multline}
with $h=h(\varkappa)$ being the characteristic determinant of problem \eqref{W0Problem}  defined in \ref{subsecScattering}.

Suppose that  $\alpha\not\in \Sigma_\Psi$. Then  $h(\varkappa)$ is different from $0$, and so
bounded as $\varepsilon\to 0$ eigenfrequencies $\omega_\varepsilon$ of \eqref{SolvModel1} tend to the roots of equation $$\tng a\omega\tng b\omega=0.$$ The corresponding eigenfunctions $y_\varepsilon$ converge in $C(a,b)$ and the set of their possible limits consists of two families of functions
$$
  y_{k,1}(x)=
    \begin{cases}
    \sin\frac{\pi k}{a}(x-a),& x\in(a,0)\\
    0,& x\in(0,b)
    \end{cases},
$$
$$
 y_{k,2}(x)=
    \begin{cases}
    0& x\in(a,0)\\
    \sin\frac{\pi k}{b}(x-b),& x\in (0,b)
    \end{cases},
$$
where $k\in \mathbb{N}$.
It is easy to check that the limit eigenfrequencies and eigenfunctions correspond to the direct sum of the second
derivative operators on intervals $(a,0)$ and $(0,b)$ subject to the Dirichlet conditions.
The sum is similar to operator
$A(\alpha,\Psi)$ for $\alpha\not\in \Sigma_\Psi$,  defined in \ref{SubsecWhatIsDeltaPrime}.

Now assume $\alpha\in \Sigma_\Psi$, and so $h(\varkappa)=0$. Remark that at the same time
$h_1(\varkappa)=\tnh^2\varkappa-1<0$.
In this case the eigenfrequencies $\omega_\varepsilon$ go to the roots of equation
 $\tng b\omega=g(\varkappa)\tng a\omega$ with
$$
   g(\varkappa)=\frac{1+\tng^2\varkappa}{1-\tnh^2\varkappa}=\frac{\csh^2\varkappa}{\cos^2\varkappa}=\theta_\Psi^2 (\alpha)\quad\text{for } \alpha\in \Sigma_\Psi\cap\mathbb R_+.
$$
The eigenfunctions $y_\varepsilon$ converge in $L_2(a,b)$ as well as uniformly on each of intervals  $(a,0)$ and $(0,b)$ towards  functions of the family
\begin{equation*}
y_k(x)=
    \begin{cases}
       r(\omega_k)\sin\omega_k(x-a),&x\in(a,0)\\
    \:\:\theta(\alpha)\sin\omega_k(x-b),&x\in(0,b)
        \end{cases},
\end{equation*}
where $\omega_k$ is a root of  $\tng b\omega=\theta_\Psi^2 (\alpha)\tng a\omega$, and
$r(\omega)=\frac{\sin b\omega}{\sin a\omega}$ for $\sin a\omega\neq 0$,
$r(\omega)=\frac{b\cos b\omega}{a\cos a\omega}$ for $\sin a\omega= 0$.
We check at once that $\omega_k^2$ and  $y_k$ are eigenvalues and eigenfunctions of the problem
\begin{equation*}
\begin{cases}
  -y''=\lambda y, \quad x\in (a,0)\cup (0,b),\qquad y(a)=y(b)=0,\\
  \phantom{-}y(+0)=\theta_\Psi (\alpha)y(-0),\quad \theta_\Psi (\alpha)y'(+0)=y'(-0),
\end{cases}
\end{equation*}
The operator associated with the problem is similar to operator
$A(\alpha,\Psi)$ in the resonant case.

\subsection{Open Problem}

In Section~\ref{SectionMainTerms} we have constructed the asymptotics in the case of an arbitrary profile belonging to $\mathcal{P}$. The case of $\delta'$-like potentials $\Psi\in \mathcal{P}_0$ is characterized by a special structure of resonance set $\Sigma_\Psi$  described in Theorem~\ref{ThSpectrumOfT}\textit{(ii)} as well as a behavior of  coupling function $\theta_\Psi$. The exactly solvable models above and the computer simulation of more complicated models suggest the following
\begin{g_hyp} Suppose that potential $\Psi$ fulfils conditions $m_0(\Psi)=0$, $m_1(\Psi)\break=-1$.
Then the coupling function $\theta_\Psi$   possesses   the following property:
\begin{itemize}\itemindent20pt
  \item[$\diamond$] $|\theta_\Psi(\alpha)|>1$ for $\alpha\in \Sigma_\Psi\cap \mathbb R_+$ and $|\theta_\Psi(\alpha)|\to +\infty$ as $\alpha\to +\infty$,
  \item[$\diamond$] $|\theta_\Psi(\alpha)|<1$ for $\alpha\in \Sigma_\Psi\cap \mathbb R_-$ and $|\theta_\Psi(\alpha)|\to 0$ as $\alpha\to -\infty$.
\end{itemize}
\end{g_hyp}
The property has its origins in the structure of eigenfunctions of
$J$-selfadjoint and $J$-nonnegative operator $\mathcal{T}_\Psi$. Because of this, the demonstration could be obtained using the Krein space theory. Note that our hypothesis for the behavior of  $\theta_\Psi$  fails for other classes of potentials. There exist even functions $\Psi$ such that $m_0(\Psi)\neq 0$ as well as
$m_0(\Psi)=m_1(\Psi)=0$. But  for an even potential we have $|\theta_\Psi(\alpha)|=1$ for all $\alpha\in \Sigma_\Psi$.

The physical interpretation of a coupling function is fairly simple.
Suppose that $v$ is  a normalized in $L_2(\mathbb{R})$ eigenfunction of  \eqref{LimitProblemV} and
$P_v(a,b)=\int_a^b |v(x)|^2\,dx$ is the probability of finding a particle in  interval $(a,b)$ provided that the system is in a pure state $v$. Then
\begin{equation*}
    \theta_\Psi ^2(\alpha)=\lim\limits_{r\to +0}\frac{P_v(0,r)}{P_v(-r,0)},
\end{equation*}
which is to say,  $\theta_\Psi ^2(\alpha)$ is the marginal ratio of probabilities of finding a particle in
$(0,r)$ and $(-r, 0)$ respectively. Let us consider potential $\frac{\alpha}{\varepsilon^{2}}\Psi(\varepsilon^{-1}x)$ with the profile  being an odd function as shown in Pic.~\ref{FigPotential} from the left.
The hypothesis states that the probability to find an electron nearby the origin is higher from the side of the deep well than that of the high wall. In other words, for  positive $\alpha$ a negative drop in the electron density occurs, while passing across the junction from the barrier to the well side. If $\alpha$ is negative, then positions of the barrier and the well are changed over.

\section{Asymptotics of Bounded Spectrum of  $\mathcal{H_{\varepsilon}}(\alpha,\Psi)$: Correctors}\label{SectionCorrector}
We need to find  next terms of  series \eqref{ExpanEVLow}--\eqref{ExpanEFWLow} in order to justify
the proximity of  energy levels of Hamiltonians $\mathcal{H_{\varepsilon}}(\alpha,\Psi)$ and $\mathcal{H}(\alpha,\Psi)$.

\subsection{Case of Closed $\delta'$-Barrier}
If coupling constant $\alpha$ doesn't belong to $\Sigma_\Psi$, then $w=0$.
Let $\lambda $ be a simple eigenvalue of  $\mathcal{H}(\alpha,\Psi)=S_-\oplus S_+$ with  $L_2(\mathbb{R})$-normalized eigenfunction $v$. Without loss of generality
we can assume that $\lambda\in \sigma(S_+)$. Clearly then $v$ vanishes in $\mathbb{R}_-$. Recalling \eqref{W1Equations} and \eqref{AdCondW1Prime} one obtains
\begin{equation*}
-w_1''+\alpha\Psi w_1=0,\quad\xi\in(-1,1),\qquad
w_1'(-1)=0,\quad w_1'(1)=v'(+0).
\end{equation*}
The problem admits a unique solution, since $\alpha$ is not an eigenvalue of \eqref{W0Problem}.
In view of \eqref{V1Equations}, \eqref{AdCondVW1}  function $v_1$
is a solution of the problems
\begin{align}\label{ProblemV1-}
  &-v_1''+U v_1=\lambda v_1, \quad x\in \mathbb{R}_-,&&v_1(-0)=w_1(-1),\\\label{ProblemV1+}
  &-v_1''+U v_1=\lambda v_1+\lambda_1 v, \quad x\in \mathbb{R}_+, &&v_1(+0)=w_1(1)-v '(+0)&
\end{align}
 on each half-axis.
The first of problems admits a unique solution in $L_2(\mathbb{R}_-)$ since  $\lambda\not\in \sigma(S_-)$.
The solution of the second one does not generally exist.
According to Fredholm's alternative problem \eqref{ProblemV1+} has a solution in $L_2(\mathbb{R}_+)$ if and only if
$\lambda_1=v '(+0)(v '(+0)-w_1(1))$. To obtain this, we  multiply  equation \eqref{ProblemV1+} by eigenfunction $v$ and then integrate twice by parts.
Solution $v_1$ is ambiguously determined, namely it is defined up to term $cv$. Subordinating it to the condition
$\int_{\mathbb{R}_-}v v_1\,dx=0$ we fix it uniquely.
 Recalling \eqref{W2Equations}, \eqref{AdCondV1W2} and $w=0$,  we obtain the problem
\begin{equation}\label{ProblemW2}
\begin{gathered}
    -w_2''+\alpha\Psi w_2=0 ,\quad \xi\in(-1,1),\\[5pt]
   w_2'(-1)=v_1'(-0),\quad w_2'(1)=v_1'(+0)+v ''(+0),
   \end{gathered}
\end{equation}
which gives us the corrector $w_2$.

We introduce the notation
\begin{equation}\label{Approx}
    \Lambda_\varepsilon=\lambda +\varepsilon\lambda_1, \quad Y_\varepsilon(x)=
\begin{cases}
v (x)+\varepsilon v_1(x),  &\abs{x}>\varepsilon,\\
 \varepsilon w_1(\varepsilon^{-1}x)+\varepsilon^2 w_2(\varepsilon^{-1}x),&\abs{x}<\varepsilon
\end{cases}
\end{equation}
for the constructed approximations of  eigenvalues and eigenfunctions.
Such approximations can also be found  in the case, when either
$\lambda\in \sigma(S_-)\setminus\sigma(S_+)$ or $\lambda\in \sigma(S_-)\cap\sigma(S_+)$.

\subsection{Case of Open $\delta'$-Barrier}
 We now assume that  $\alpha$ belongs to  resonant set $\Sigma_\Psi$, and $w= a W_\alpha$, where  $W_\alpha$ is an eigenfunction of  \eqref{W0Problem} and $a$ is an arbitrary constant.
In this case $\lambda$ is an eigenvalue of \eqref{LimitProblemV} with
$L_2(\mathbb{R})$-normalized eigenfunction $v$.
Note that all eigenvalues of this problem are simple.
The first condition \eqref{AdCondVW} yields $a=\frac{v(-0)}{W_\alpha(-1)}$.
Then the second one also holds since
$ v(+0)=\theta_\Psi(\alpha)v(-0)=\frac{v(-0)}{W_\alpha(-1)}W_\alpha(1)=aW_\alpha(1)$.
In view of Remark \ref{RemarkV(0)=0}  constant $a$ is zero if and only if $\lambda\in \sigma(S_-)\cap\sigma(S_+)$.

Problem \eqref{W1Problem} admits  solution $w_1$, because \eqref{ThetaVPrimeCond} holds.
This solution can be represented as $w_1=w_1^*+a_1W_\alpha$, where $w_1^*$ is a partial solution of the problem and $a_1$ is an arbitrary constant.
This constant will be found below, but first we shall construct  $v_1$.

Function $v_1$ satisfies  equation \eqref{V1Equations} outside of the origin and
\begin{equation*}
    v_1(+0)-\theta_\Psi(\alpha)v_1(-0)=g_1
\end{equation*}
with $g_1=w_1(1)-\theta_\Psi(\alpha)w_1(-1)-v '(+0)-\theta_\Psi(\alpha)v '(-0)$, which follows from equalities \eqref{AdCondVW1}.
Although $a_1$ is as yet unknown, constant $g_1$ is uniquely determined. Indeed,
\begin{multline*}
    w_1(1)-\theta_\Psi(\alpha)w_1(-1)=w_1^*(1)-\theta_\Psi(\alpha)w_1^*(-1)+\\
    +a_1\bigl(W_\alpha(1)-\theta_\Psi(\alpha)W_\alpha(-1)\bigr)=w_1^*(1)-\theta_\Psi(\alpha)w_1^*(-1).
\end{multline*}
Next we employ \eqref{W2Equations} and \eqref{AdCondV1W2} to derive the problem
\begin{equation}\label{W2Problem}
\begin{cases}
-w_2''+\alpha\Psi(\xi)w_2=(\lambda-U(0))w ,\qquad \xi\in(-1,1),\\
\phantom{-}w_2'(-1)=v_1'(-0)-v ''(-0),\quad w_2'(1)=v_1'(+0)+v ''(+0).
\end{cases}
\end{equation}
Its solvability condition can be written as
\begin{equation}\label{SolvabilityW2}
    \theta_\Psi(\alpha)v_1'(+0)-v_1'(-0)=h_1
\end{equation}
with $h_1=(w(-1))^{-1}(\lambda -U(0))\int_{-1}^1w^2\,d\xi -\theta_\Psi(\alpha)v ''(+0)-v ''(-0)$.
Therefore  $v_1$ must be a solution of the problem
\begin{equation}\label{ProblemV2}
\begin{cases}
    -v_1''+Uv_1=\lambda v_1+\lambda_1v,\qquad x\in\mathbb R\setminus\{0\},\\
    \phantom{-}v_1(+0)-\theta_\Psi(\alpha)v_1(-0)=g_1,\quad \theta_\Psi(\alpha)v_1'(+0)-v_1'(-0)=h_1.
\end{cases}
\end{equation}
Free parameter $\lambda_1$ in the right-hand side of equation \eqref{ProblemV2} enables us to solve the problem.
In light of Fredholm's alternative, \eqref{ProblemV2} admits a solution if and only if
$\lambda_1=g_1v '(-0)-h_1v (-0)$. For definiteness sake, the solution is subject to additional condition
$\int_{\mathbb{R}}v v_1\,dx=0$.

Given $v_1$, we can now compute constant $a_1$. From the first condition in \eqref{AdCondVW1} we deduce  $a_1=(W_\alpha(-1))^{-1}\bigl(v_1(-0)-v'(-0)-w_1^*(-1)\bigr)$. It is easy to check directly that the second condition in \eqref{AdCondVW1} also holds.

Finally in the case of  open $\delta'$-barrier one obtains such approximations for the eigenvalue and the eigenfunction of the perturbed problem
\begin{equation}\label{Approx1}
\begin{gathered}
    \Lambda_\varepsilon=\lambda +\varepsilon\lambda_1, \\ Y_\varepsilon(x)=
\begin{cases}
v (x)+\varepsilon v_1(x),  &\abs{x}>\varepsilon,\\
\frac{v(-0)}{W_\alpha(-1)}W_\alpha(\varepsilon^{-1}x) + \varepsilon
w_1(\varepsilon^{-1}x)+\varepsilon^2 w_2(\varepsilon^{-1}x),&\abs{x}<\varepsilon.
\end{cases}
\end{gathered}
\end{equation}
Here $w_2$ is an arbitrary solution of  \eqref{W2Problem}. Recall that condition \eqref{SolvabilityW2}
ensures its existence.
The choice of  $a_2$ in  representation $w_2=w_2^*+a_2W_\alpha$ is of no importance since we do not look for corrector $v_2$.

\section{Justification of Asymptotic Expansions}\label{SectionMainResult}
As shown in  Theorem~\ref{TheoremNeg},
for some singular potential $\alpha \Psi_\varepsilon(x)$ there is  finite number
$N(\alpha,\Psi)>0$ of eigenvalues $\lambda_k^\varepsilon(\alpha, \Psi)$, converging to $-\infty$ as $\varepsilon\to 0$. Other eigenvalues remain bounded as $\varepsilon\to 0$.  We shall show that these eigenvalues converge to the eigenvalues of  $\mathcal{H}(\alpha,\Psi)$.

\subsection{Convergence Theorem}
Let $\{\lambda_\varepsilon\}_{\varepsilon\in \mathcal{I}}$ be a sequence of eigenvalues of  $\mathcal{H}_{\varepsilon}(\alpha,\Psi)$ and  assume that $\{y_\varepsilon\}_{\varepsilon\in \mathcal{I}}$ is a sequence of the corresponding $L_2(\mathbb{R})$-normalized  eigenfunctions. Here $\mathcal{I}$ is an infinite subset of $(0,1)$ for which $0$ is an accumulation point.

\begin{g_theorem}\label{ConvergenceTheorem}
If $\lambda_\varepsilon\to \lambda$ and $y_\varepsilon \to v$ in $L_2(\mathbb{R})$ weakly as $\mathcal{I}\ni \varepsilon\to 0$, then $\lambda$ is an eigenvalue of $\mathcal{H}(\alpha,\Psi)$ with the corresponding eigenfunction $v$. Furthermore, $y_\varepsilon$ converges to $v$ in the $L_2(\mathbb{R})$-norm.
\end{g_theorem}

We have divided the proof into a sequence of lemmas.

\begin{g_lemma}\label{LemConvergenceC1}
Under the assumptions of Theorem~\ref{ConvergenceTheorem}, for every compact $K\subset \mathbb{R}$ and $\gamma>0$ the sequence $y_\varepsilon$ converges to $v$ in $W_2^2(K\setminus(-\gamma,\gamma))$ weakly as well as in the norm of $C^1(K\setminus(-\gamma,\gamma))$.
Moreover, v is a solution of the Schr\"{o}dinger
equation
\begin{equation}\label{UnperturbedSch}
-v''+Uv=\lambda v
\end{equation}
on each half-axis.
\end{g_lemma}
\begin{proof}
Let $\mathcal{M}_\gamma$ be the set of test functions  $\varphi\in C_0^\infty(\mathbb{R})$ such that $\varphi(x)=0$ for $x\in (-\gamma/2,\gamma/2)$.
Since the function $y_\varepsilon$ is smooth, it belongs to the space $W_{2,loc}^2(\mathbb{R})$.
We conclude from \eqref{MainProbl} that
\begin{equation}\label{IntegIdentEps}
    \int_\mathbb{R} y_\varepsilon''\varphi\,dx=\int_\mathbb{R} (U-\lambda_\varepsilon)y_\varepsilon\varphi\,dx,
\end{equation}
for all  $\varphi\in \mathcal{M}_\gamma$ and $\varepsilon<\gamma/2$, since the support of the short-range potential $\Psi_\varepsilon$ lies in $(-\gamma/2,\gamma/2)$.
The right-hand side of \eqref{IntegIdentEps} has a limit as $\varepsilon\to 0$ by the assumptions.
Thus, the left-hand side also converges for all  $\varphi\in \mathcal{M}_\gamma$.
We deduce then that $y_\varepsilon \to v$ in $W_{2, loc}^2(\mathbb{R}\setminus[-\gamma/2,\gamma/2])$ weakly, and so
in $W_2^2(K\setminus(-\gamma,\gamma))$ weakly.
Moreover,
\begin{equation*}\label{IntegIdent}
    \int_\mathbb{R} v''\varphi\,dx=\int_\mathbb{R} (U-\lambda)v\varphi\,dx \qquad\text{for all}\quad \varphi\in \mathcal{M}_\gamma.
\end{equation*}
It follows hereby that $v$ is a solution of \eqref{UnperturbedSch} on $\mathbb{R}\setminus(-\gamma,\gamma)$ and, therefore,
on whole half-lines in view of the arbitrariness of $\gamma$.
It remains to note that the imbedding theorem implies the convergence of $y_\varepsilon$ in $C^1(K\setminus(-\gamma,\gamma))$.
\end{proof}

\begin{g_lemma}\label{LemmmaYe(eps)}
Suppose as in Theorem~~\ref{ConvergenceTheorem} that $\lambda_\varepsilon\to \lambda$ and $y_\varepsilon \to v$ in $L_2(\mathbb{R})$ weakly as $\mathcal{I}\ni \varepsilon\to 0$.
Then $y_\varepsilon(\varepsilon)\to v(+0)$, $y_\varepsilon'(\varepsilon)\to v'(+0)$, $y_\varepsilon(-\varepsilon)\to v(-0)$ and $y_\varepsilon'(-\varepsilon)\to v'(-0)$ as $\mathcal{I}\ni \varepsilon\to 0$.
\end{g_lemma}
\begin{proof}
Let $\zeta$ be the $C^\infty(\mathbb{R}\setminus \{0\})$-function such that $\zeta(x)=0$ for $x<0$ and $x> 2$, and also $\zeta(x)=1$ for $x\in (0,1)$. We hereafter denote the characteristic function of set $K$ by $\chi_K$. Let us introduce  sequence   $\zeta_\varepsilon(x)=\chi_{(\varepsilon,\infty)}(x)\zeta(x)$. Multiplying both equations \eqref{MainProbl}, \eqref{UnperturbedSch} by $\zeta_\varepsilon$ and integrating by parts yield
$$
y_\varepsilon'(\varepsilon)=-\int\limits^{\infty}_1y_\varepsilon'\zeta'\,dx+
\int\limits^{\infty}_\varepsilon(\lambda_\varepsilon-U)y_\varepsilon\zeta\,dx,
$$
$$
v'(\varepsilon)=-\int\limits^{\infty}_1v'\zeta'\,dx+
\int\limits^{\infty}_\varepsilon(\lambda-U)v\zeta\,dx.
$$
Here we take into account that $\zeta_\varepsilon(\varepsilon+0)=1$ and $\zeta'_\varepsilon(x)=0$ for $x\in (\varepsilon,1)$. According to Lemma~\ref{LemConvergenceC1}
the right-hand sides of the equalities have the same limit as $\varepsilon\to 0$, then $y_\varepsilon'(\varepsilon)\to v'(+0)$.
In a similar way, using  cut-function $\zeta_\varepsilon(-x)$, we can prove that $y_\varepsilon'(-\varepsilon)\to v'(-0)$.

Now let us consider the cut-function $\eta\in C^\infty(\mathbb{R}\setminus \{0\})$ such that $\eta(x)=0$ for $x<0$ and $x>2$, and also $\eta(x)=x$ for $x\in (0,1)$. We set $\eta_\varepsilon(x)=\chi_{(\varepsilon,\infty)}(x)\eta(x)$. In view of equalities $\eta_\varepsilon(\varepsilon+0)=\varepsilon$ and $\zeta''_\varepsilon(x)=0$ for $x\in (\varepsilon,1)$, and employing  \eqref{MainProbl}, \eqref{UnperturbedSch} we derive
\begin{gather*}
y_\varepsilon(\varepsilon)=\varepsilon y_\varepsilon'(\varepsilon)
-\int\limits^{\infty}_1y_\varepsilon\eta''\,dx+
\int\limits^{\infty}_\varepsilon(U-\lambda_\varepsilon)y_\varepsilon\eta\,dx,\\[5pt]
v(\varepsilon)=\varepsilon v'(\varepsilon)-\int\limits^{\infty}_1v\eta''\,dx+
\int\limits^{\infty}_\varepsilon(U-\lambda)v\eta\,dx.
\end{gather*}
Because $y_\varepsilon'(\varepsilon)\to v'(+0)$, we conclude that  the right-hand sides of these equalities converge to the same limit, and so $y_\varepsilon(\varepsilon)\to v(+0)$.
The same proof works for the sequence $y_\varepsilon(-\varepsilon)$.
\end{proof}

We are now in a position to describe the behavior of  eigenfunctions $y_\varepsilon$ in a neighborhood of the origin. Let $w$ and $z$ be solutions of the Cauchy problems on $[-1,1]$
\begin{align}\label{CauchyProblemW}
   & -w''+\alpha\Psi(\xi) w=0,\qquad w(-1)=1, \quad  w'(-1)=0;\\ \label{CauchyProblemZ}
   & -z''+\alpha\Psi(\xi) z=0,\qquad     z(-1)=0,\quad     z'(-1)=v'(-0).
\end{align}

\begin{g_lemma}\label{LemmaLocalConvergence} If $\lambda_\varepsilon\to \lambda$ and $y_\varepsilon \to v$ in $L_2(\mathbb{R})$ weakly as $\mathcal{I}\ni \varepsilon\to 0$, then
  \begin{equation}\label{YeWZ}
  \left\|\varepsilon^{-1}y_\varepsilon(\varepsilon\xi)-\varepsilon^{-1}y_\varepsilon(-\varepsilon)w(\xi)- z(\xi)\right\|_{C^1([-1,1])}\to 0.
  \end{equation}
\end{g_lemma}
\begin{proof}
Set $w_\varepsilon(\xi)=\varepsilon^{-1}y_\varepsilon(\varepsilon\xi)
-\varepsilon^{-1}y_\varepsilon(-\varepsilon)w(\xi)- z(\xi)$.
Upon substituting $\xi=\varepsilon^{-1}x$ we see that equation \eqref{MainProbl} may be rewritten as
\begin{equation*}
    -\frac{d^2y_\varepsilon}{d\xi^2}+\alpha\Psi(\xi)y_\varepsilon=\varepsilon^2 (\lambda_\varepsilon-U(\varepsilon\xi))y_\varepsilon.
  \end{equation*}
Employing \eqref{CauchyProblemW} and \eqref{CauchyProblemZ} we deduce that $w_\varepsilon$ is a solution of the Cauchy problem
  \begin{equation}\label{CauchyProblemWe}
    \begin{cases}
    -w''_\varepsilon+\alpha\Psi(\xi) w_\varepsilon=f_\varepsilon(\xi),\qquad \xi\in[-1,1],\\
    w_\varepsilon(-1)=0,\quad     w'_\varepsilon(-1)=y_\varepsilon'(-\varepsilon)-v'(-0)
    \end{cases}
\end{equation}
with $f_\varepsilon(\xi)=\varepsilon(\lambda_\varepsilon-U(\varepsilon\xi))y_\varepsilon(\varepsilon\xi)$. For every $f_\varepsilon\in L_2(-1,1)$ there exists a unique solution $w_\varepsilon\in W_2^2(-1,1)$. Moreover this solution satisfies the estimate
\begin{equation}\label{AprioriEst}
    \|w_\varepsilon\|_{W_2^2(-1,1)}\leq C\left(\|f_\varepsilon\|_{L_2(-1,1)}+|y_\varepsilon'(-\varepsilon)-v'(-0)|\right)
\end{equation}
with constant $C$ being independent of $\varepsilon$. Next the inequality
\begin{equation*}
    \int_{-1}^1y_\varepsilon^2(\varepsilon\xi)\,d\xi=\frac{1}{\varepsilon}\int_{-\varepsilon}^\varepsilon y_\varepsilon^2(x)\,dx\leq \frac{1}{\varepsilon}\,\|y_\varepsilon\|_{L_2(\mathbb{R})}=\frac{1}{\varepsilon}
\end{equation*}
implies $\varepsilon^{-1/2}\|y_\varepsilon(\varepsilon\xi)\|_{L_2(-1,1)}\leq c$. Hence $\|f_\varepsilon\|_{L_2(-1,1)}\leq c \varepsilon^{1/2}$, because $\lambda_\varepsilon$ is bounded.
Therefore in light of Lemma~\ref{LemmmaYe(eps)} the right-hand side of \eqref{AprioriEst} tends to zero.
The proof is complete by using the imbedding theorem $W_2^2(-1,1)\subset C^1([-1,1])$.
\end{proof}

\begin{g_lemma}\label{LemmaL2convergence}
If $\lambda_\varepsilon\to\lambda$ and $y_\varepsilon\to v$ in $L_2(\mathbb{R})$ weakly as $\mathcal{I}\ni\varepsilon\to0$, then $y_\varepsilon\to v$ in $L_2(\mathbb{R})$.
\end{g_lemma}
\begin{proof}
Let us first show that $y_\varepsilon$ is uniformly bounded on each compact $K=[-A,A]$, $A>0$.
Applying \eqref{YeWZ} and Lemma~\ref{LemmmaYe(eps)} we see at once that the sequence $y_\varepsilon$
is uniformly bounded on $[-\varepsilon,\varepsilon]$.
Indeed,
\begin{equation}\label{Uniform-ee}
\max\limits_{x\in[-\varepsilon,\varepsilon]}|y_\varepsilon(x)|\leq c\varepsilon+|z(\varepsilon^{-1}x)|\varepsilon+|y_\varepsilon(-\varepsilon)|\,|w(\varepsilon^{-1}x)|\leq c_1.
\end{equation}
Let $\Omega_\varepsilon=K\setminus (-\varepsilon,\varepsilon)$.
Multiplying equation~\eqref{MainProbl} by $\chi_{\Omega_\varepsilon}y_\varepsilon$ and integrating
by parts give us
\begin{equation*}
   \int_{\Omega_\varepsilon}y_\varepsilon'^{2}\,dx=
\int_{\Omega_\varepsilon}(\lambda_\varepsilon-U)y_\varepsilon^2\,dx+y_\varepsilon'y_\varepsilon\Big|_{-A}^{-\varepsilon}+y_\varepsilon'y_\varepsilon\Big|^{A}_{\varepsilon}.
\end{equation*}
According to Lemmas \ref{LemConvergenceC1} and \ref{LemmmaYe(eps)}, all terms on the right-hand side are uniformly bounded with respect to $\varepsilon$.
From this we conclude that the sequence $y_\varepsilon$ is bounded in $W^1_2(\Omega_\varepsilon)$ by a constant being independent of $\varepsilon$, hence that  $\max\limits_{x\in\Omega_\varepsilon} |y_\varepsilon(x)|\leq c_2$, and finally that $\max\limits_{x\in K}|y_\varepsilon(x)|\leq c_3$, by
\eqref{Uniform-ee}.
Here the constant $c_3$ does not depend on $\varepsilon$.

The potential $U$ increases for $x\to\pm\infty$, then for large $\abs{x}$ the eigenfunctions $y_\varepsilon$ and $v$ are exponentially small.
In particular,
$$
    |y_\varepsilon(x)|\leq e^{-|x|},\qquad|v(x)|\leq e^{-|x|}
$$
on each set $|x|>A$, where the potential $q-\lambda_\varepsilon$ and $q-\lambda$ are greater than one \cite[p.~59]{BS}.

Fix $\delta>0$ and choose a constant $A$ so large that  $\|y_\varepsilon\|_{L_2(|x|>A)}<\delta$ and $\|v\|_{L_2(|x|>A)}<\delta$.
We see from Lemma~\ref{LemConvergenceC1} that for small enough $\varepsilon$ the $L_2(\gamma<|x|<A)$-norm
of the difference $y_\varepsilon-v$ is smaller than $\delta$.
Moreover,
\begin{equation*}
    \|y_\varepsilon-v\|_{L_2(|x|<\gamma)}\leq \sqrt{2\gamma}\,\max_{x\in [-\gamma,\gamma] }\left(|y_\varepsilon(x)|+|v(x)|\right)\leq c_4\sqrt{\gamma}<\delta
\end{equation*}
for small $\gamma$.
Thus,
\begin{multline*}
   \|y_\varepsilon-v\|_{L_2(\mathbb{R})}\leq  \|y_\varepsilon\|_{L_2(|x|>A)}+\|v\|_{L_2(|x|>A)}+\\+\|y_\varepsilon-v\|_{L_2(\gamma<|x|<A)}
   + \|y_\varepsilon-v\|_{L_2(|x|<\gamma)}< 4\delta
\end{multline*}
for small enough $\varepsilon$, $\gamma$ such that $\varepsilon<\gamma$.
It remains to note that $\delta$ may be chosen arbitrary small.
\end{proof}

\begin{proof}[Proof of  Theorem $\ref{ConvergenceTheorem}$.]
We conclude from Lemmas~\ref{LemConvergenceC1}, \ref{LemmaL2convergence} that  $v$ is a solution of the equation
\begin{equation*}
    -v''+Uv=\lambda v,\qquad x\in \mathbb{R}\setminus\{0\}
\end{equation*}
and $\|v\|_{L_2(\mathbb R)}=1$. We are left with the task of showing that $v$ satisfies appropriate coupling conditions at the origin.
In light of lemma~\ref{LemmaLocalConvergence}, we thus deduce
\begin{equation}\label{PerelimitConds}
    \varepsilon^{-1}\bigl(y_\varepsilon(\varepsilon)-y_\varepsilon(-\varepsilon)w(1)\bigl)\to z(1),\quad
    y_\varepsilon'(\varepsilon)-\varepsilon^{-1}y_\varepsilon(-\varepsilon)w'(1)\to z'(1)
\end{equation}
as $\varepsilon\to 0$. According to Lemma~\ref{LemmmaYe(eps)}, we have $y_\varepsilon(\varepsilon)-y_\varepsilon(-\varepsilon)w(1)\to 0$, $y_\varepsilon(-\varepsilon)w'(1)\to 0$. Therefore
\begin{gather}\label{LimitCondForV0}
    v(+0)-v(-0)w(1)=0,\\\label{LimitCondForV1}
    v(-0)w'(1)=0.
\end{gather}

First let assume that $v(-0)\neq 0$. Then $w'(1)=0$. By \eqref{CauchyProblemW}  $w$ is a nontrivial solution of the boundary value problem $-w''+\alpha\Psi w=0$, $w'(-1)=0$, $w'(1)=0$. Hence $w$ is an eigenfunction of operator $\mathcal{T}_\Psi$ and coupling constant $\alpha$ belongs to  resonant set $\Sigma_\Psi$. Furthermore $\theta_\Psi (\alpha)=w(1)$ since $w(-1)=1$ by construction. Thus \eqref{LimitCondForV0} coincides with coupling condition \eqref{ThetaVCond}. We obtain from the second formula in  \eqref{PerelimitConds} that $v'(+0)=z'(1)$, and so $z$ is a solution of the boundary value problem
\begin{equation}\label{ZBoundaryValueProblem}
    -z''+\alpha\Psi z=0,\qquad z'(-1)=v'(-0),\quad z'(1)=v'(+0),
\end{equation}
which coincides with \eqref{W1Problem}.
Then  solvability condition  \eqref{ThetaVPrimeCond} must hold, which is the second coupling condition for $v$.
Thus $v$ is an eigenfunction of $\mathcal{H}(\alpha,\Psi)$ with the eigenvalue $\lambda$. Here $\lambda$ is a limit of sequence $\lambda_\varepsilon$.

Now assume $v(-0)=0$. Employing \eqref{LimitCondForV0} yields $v(+0)=0$. If $\alpha\not\in \Sigma_\Psi$, then the proof is complete. In the case, when
$\alpha$ belongs to $\Sigma_\Psi$, function $v$ satisfies  condition \eqref{ThetaVCond}  and $w$  is an eigenfunction of  $\mathcal{T}_\Psi$. Then there exists  solution $z$  of \eqref{ZBoundaryValueProblem}, consequently
 condition \eqref{ThetaVPrimeCond} holds.
This case is described in Remark~\ref{RemarkV(0)=0}.
\end{proof}

\begin{g_corollary}\label{Corr1}
Suppose that  eigenvalue $\lambda_k^\varepsilon(\alpha, \Psi)$ of $\mathcal{H}_{\varepsilon}(\alpha,\Psi)$ is bounded from below. Then $\lambda_k^\varepsilon(\alpha, \Psi)$ has a finite limit as $\mathcal{I}\ni \varepsilon\to 0$
and this limit is a point of the spectrum  of $\mathcal{H}(\alpha,\Psi)$.
\end{g_corollary}
\begin{proof}
Suppose, contrary to our claim, that
\begin{equation*}
    \mu_*=\varliminf_{\varepsilon\to 0}\lambda_k^\varepsilon(\alpha, \Psi)<\varlimsup_{\varepsilon\to 0}\lambda_k^\varepsilon(\alpha, \Psi)=\mu^*.
  \end{equation*}
The constants $\mu_*$, $\mu^*$ are finite since $\lambda_k^\varepsilon(\alpha, \Psi)$ is a bounded function.
Recall that $\lambda_k^\varepsilon(\alpha, \Psi)$ is a continuous function of $\varepsilon\in (0,1)$.
Then for each $\lambda\in [\mu_*,\mu^*]$ there exists  a subsequence of eigenvalues $\lambda_\varepsilon$, $\varepsilon\in \mathcal{I}$, converging to $\lambda$. For instance, set $\mathcal{I}$ can be chosen
as a sequence of the roots of  equation  $\lambda_k^\varepsilon(\alpha, \Psi)=\lambda$ with respect to $\varepsilon$.
Sequence $\{y_\varepsilon\}_{\varepsilon\in \mathcal{I}}$ of the corresponding normalized eigenfunctions contains a weakly convergent subsequence. By the theorem  $\lambda$ is an eigenvalue of  $\mathcal{H}(\alpha,\Psi)$.
Therefore interval $[\mu_*,\mu^*]$ belongs to spectrum $\sigma(\mathcal{H}(\alpha,\Psi))$, but this is impossible unless $\mu_*=\mu^*$.
\end{proof}

\begin{g_corollary}\label{Corr2} For each eigenvalue $\lambda$ of  $\mathcal{H}(\alpha,\Psi)$ with multiplicity $s$ there exist exactly $s$ eigenvalues $\lambda_k^\varepsilon(\alpha, \Psi)$ of $\mathcal{H}_{\varepsilon}(\alpha,\Psi)$
converging to $\lambda$ as $\varepsilon\to 0$.
\end{g_corollary}
\begin{proof} Operator $\mathcal{H}(\alpha,\Psi)$ has  simple or double eigenvalues.
Let us assume  that $s=1$, but  $\lambda_k^\varepsilon(\alpha, \Psi)\to \lambda$ and  $\lambda_{k+1}^\varepsilon(\alpha, \Psi)\to \lambda$ for some $k$. Then there exist two sequences $\{y_k^{\varepsilon}(x;\alpha, \Psi)\}_{\varepsilon\in \mathcal{I}}$, $\{y_{k+1}^{\varepsilon}(x;\alpha, \Psi)\}_{\varepsilon\in \mathcal{I}}$ of the eigenfunctions, which  converge in $L_2(\mathbb{R})$
to   vectors  of the form $e^{i\varphi} v$. This contradicts the fact that $y_k^{\varepsilon}(x;\alpha, \Psi)$
and $y_{k+1}^{\varepsilon}(x;\alpha, \Psi)$ are orthogonal in $L_2(\mathbb{R})$ for all $\varepsilon\in \mathcal{I}$. Similar arguments apply to the case $s=2$ (see \cite{GolovatyMMO} for more details).
\end{proof}

\subsection{Quasimodes of $\mathcal{H}_{\varepsilon}(\alpha,\Psi)$}
We now show that each  point of $\sigma(\mathcal{H}(\alpha,\Psi))$ is a limit of
the eigenvalues of  $\mathcal{H}_{\varepsilon}(\alpha,\Psi)$.
Let $B$ be a self-adjoint operator in Hilbert space  $H$ with a domain
$\mathcal{D}(B)$. A pair $(\mu, u)\in \mathbb{R}\times\mathcal{D}(B)$ with $\|u\|_H=1$
is called a \emph{quasimode} of the operator $B$ with an accuracy up to $\rho>0$ if $\|Bu-\mu u\|_H\leq \rho$.

\begin{g_lemma}\label{LemmaVishik} Suppose that the spectrum of $B$ is discrete and simple. If $(\mu, u)$ is a quasimode of $B$ with accuracy to $\rho>0$, then interval  $[\mu-\rho,\mu+\rho]$ contains an eigenvalue $\lambda$ of $B$.
Furthermore, if segment $[\mu-\tau,\mu+\tau]$ contains only this eigenvalue of $B$, then $\|u-v\|_H\leq 2\tau^{-1}\rho$, where $v$ is an eigenfunction of $B$ for the eigenvalue $\lambda$ and $\|v\|_H=1$.
\rm{\cite{LazutkinVINITI34, VishykLust}}
\end{g_lemma}

Inasmuch as we didn't construct the asymptotics for the case of double eigenvalue $\lambda$ in the previous section, we assume that $\sigma(S_-)\cap\sigma(S_+)=\emptyset$ which ensures that  the spectrum of $\mathcal{H}(\alpha,\Psi)$ is simple.

Let us construct the quasimodes of $\mathcal{H}_{\varepsilon}(\alpha,\Psi)$.
Now fix some eigenvalue $\lambda$ of  $\mathcal{H}(\alpha,\Psi)$ with eigenfunction $v$, $\|v\|_{L_2(\mathbb{R})}=1$. For each $\lambda$ and $v$ we have obtained the formal asymptotic approximations
$\Lambda_\varepsilon$, $Y_\varepsilon$ defined by either \eqref{Approx} or \eqref{Approx1} depending upon  $\alpha$ and $\Psi$.
The task is  to justify the asymptotics.
Notice that in this subsection we don't distinguish the resonant and non-resonant  cases.
By construction $\Lambda_\varepsilon$ and $Y_\varepsilon$ satisfy equalities
\begin{align}
\begin{aligned}\label{Reminders}
    -&Y_\varepsilon''+(U(x)-\Lambda_\varepsilon) Y_\varepsilon=\varepsilon^2 R_1(\varepsilon,x)\quad \text{for } \abs{x}>\varepsilon,\\
    -&Y_\varepsilon''+(U(x)+\alpha \varepsilon^{-2}\Psi(\varepsilon^{-1}x)-\Lambda_\varepsilon) Y_\varepsilon=\varepsilon R_2(\varepsilon,x)\quad \text{for }\abs{x}<\varepsilon,\\
    [&Y_\varepsilon]_{x=\pm\varepsilon}=\varepsilon^2 r_1^\pm(\varepsilon),\qquad [Y'_\varepsilon]_{x=\pm\varepsilon}=\varepsilon^2 r_2^\pm(\varepsilon)
\end{aligned}
\end{align}
with  functions $R_j$, $r_j^\pm$ being uniformly bounded with respect to its arguments.
Function $Y_\varepsilon$ doesn't belong to the domain of $\mathcal{H}_{\varepsilon}(\alpha,\Psi)$, because it has jump discontinuities  at points $x=\pm\varepsilon$.  But we can construct function $\zeta_\varepsilon$ with the following properties
\begin{itemize}
  \item $\zeta_\varepsilon$ is a smooth function out of points $x=\pm\varepsilon$ and differs from zero only for $\varepsilon<|x|<1$;
  \item $[\zeta_\varepsilon]_{x=\pm\varepsilon}=-r_1^\pm(\varepsilon)$ and $[\zeta'_\varepsilon]_{x=\pm\varepsilon}=-r_2^\pm(\varepsilon)$;
  \item $\max\limits_{\varepsilon<|x|<1}(|\zeta_\varepsilon(x)|+|\zeta'_\varepsilon(x)|+|\zeta''_\varepsilon(x)|)\leq c$ with constant $c$ being independent of $\varepsilon$,
\end{itemize}
which eliminates this defect. In fact, function $Y_\varepsilon+\varepsilon^2\zeta_\varepsilon$ is continuously differentiable at $x=\pm\varepsilon$ and belongs to $\mathcal{D}(\mathcal{H}_{\varepsilon}(\alpha,\Psi))$. We set $\Upsilon_\varepsilon=\|Y_\varepsilon+
\varepsilon^2\zeta_\varepsilon\|^{-1}_{L_2(\mathbb{R})}(Y_\varepsilon+\varepsilon^2\zeta_\varepsilon)$
and substitute $\Upsilon_\varepsilon$ into \eqref{Reminders} instead of $Y_\varepsilon$. Then
the orders of smallness of  right-hand sides in \eqref{Reminders} don't  change since $\|Y_\varepsilon+
\varepsilon^2\zeta_\varepsilon\|_{L_2(\mathbb{R})}\to 1$ as $\varepsilon\to 0$. Therefore  pair $(\Lambda_\varepsilon, \Upsilon_\varepsilon)$ is a quasimode of  $\mathcal{H}_{\varepsilon}(\alpha,\Psi)$ with accuracy to $\varepsilon$.

\begin{g_lemma}\label{LemmaQuasimodes}
For each $\lambda\in \sigma(\mathcal{H}(\alpha,\Psi))$ there is an eigenvalue $\lambda_j^\varepsilon(\alpha, \Psi)$ of $\mathcal{H}_{\varepsilon}(\alpha,\Psi)$ such that $\lambda_j^\varepsilon(\alpha, \Psi)\to \lambda$. Moreover,
  \begin{equation}\label{LemmaEstVishik}
 \abs{\lambda_j^\varepsilon(\alpha, \Psi)-\lambda}\leq c_1\varepsilon,\qquad
 \bigl\|y_j^{\varepsilon}(\,\cdot\,;\alpha, \Psi)-v\bigr\|_{L_2(\mathbb{R})}\leq c_2\varepsilon,
\end{equation}
where $y_j^{\varepsilon}$ and $v$ are the corresponding normalized eigenfunctions.
\end{g_lemma}
\begin{proof} Let $(\Lambda_\varepsilon, \Upsilon_\varepsilon)$ be the quasimode of $\mathcal{H}_{\varepsilon}(\alpha,\Psi)$ corresponding to the limit eigenvalue $\lambda$ and eigenfunction $v$. According to Lemma~\ref{LemmaVishik} there exists number $j$ such that $\abs{\lambda_j^\varepsilon(\alpha, \Psi)-\Lambda_\varepsilon}\leq c_1\varepsilon$, from which  the first inequality in \eqref{LemmaEstVishik} follows.
In view of Corollaries~\ref{Corr1}, \ref{Corr2}  index
$j$ is independent of $\varepsilon$. If $\tau$ is less than the distance from $\lambda$ to the rest of the spectrum of $\mathcal{H}(\alpha,\Psi)$, then  for small enough $\varepsilon$ interval $[\lambda-\tau,\lambda+\tau]$  contains eigenvalue $\lambda_j^\varepsilon(\alpha, \Psi)$ only. Applying again Lemma~\ref{LemmaVishik} yields
  $\bigl\|y_j^{\varepsilon}(\,\cdot\,;\alpha, \Psi)-\Upsilon_\varepsilon\bigr\|_{L_2(\mathbb{R})}\leq 2\tau^{-1}c_1\varepsilon$, from which  the second inequality in \eqref{LemmaEstVishik} immediately follows.
\end{proof}

We are now in a position to state the main result.
Let $\{\lambda_k(\alpha, \Psi)\}_{k=1}^\infty$ be the eigenvalues of  $\mathcal{H}(\alpha,\Psi)$ numbered by increasing  and let $\{v_k(x;\alpha, \Psi)\}_{k=1}^\infty$ be the orthonormal system of eigenfunctions. Recall that $N=N(\alpha,\Psi)$ is the number of eigenvalues of $\mathcal{H_{\varepsilon}}(\alpha,\Psi)$  converging to $-\infty$ as $\varepsilon\to 0$.
\begin{g_theorem}
Assume $(\alpha,\Psi)\in \mathbb{R}\times \mathcal{P}$ and $\sigma(S_-)\cap\sigma(S_+)=\emptyset$. Then for every natural $k$ we have
\begin{align}\label{EstLambda}
 &\abs{\lambda_{k+N}^\varepsilon(\alpha, \Psi)-\lambda_k(\alpha, \Psi)}\leq c_1\varepsilon,\\\label{EstY}
 &\bigl\|y_{k+N}^{\varepsilon}(\,\cdot\,;\alpha, \Psi)-v_k(\,\cdot\,;\alpha, \Psi)\bigr\|_{L_2(\mathbb{R})}\leq c_2\sqrt{\varepsilon}
\end{align}
with  constants $c_1$, $c_2$ being independent of $\varepsilon$.
\end{g_theorem}
\begin{proof}
Theorem~\ref{ConvergenceTheorem}, Corollary~\ref{Corr2} and Lemma~\ref{LemmaQuasimodes}
may be summarized by saying that eigenvalue
$\lambda_{k+N}^\varepsilon(\alpha, \Psi)$ tends to $\lambda_k(\alpha, \Psi)$ as $\varepsilon\to 0$.
Thus  estimates \eqref{EstLambda}, \eqref{EstY} are more precise versions of \eqref{LemmaEstVishik}.
\end{proof}

\paragraph{Acknowledgements.}
We wish to thank  Branko \'{C}urgus, Illia Kamotski and Rostyslav Hryniv for useful discussions.
We  are also greatly indebted to Oleksiy Kostenko, the referee, and most especially to Mark Malamud,
for carefully reading the paper and suggesting some improvements.

%\makepaperend

\end{document}